\documentclass[final,3p,times]{elsarticle}

\usepackage{framed,multirow}

\usepackage{amssymb}
\usepackage{amsmath}
\usepackage{latexsym}

\usepackage{url}
\usepackage[dvipsnames]{xcolor}
\definecolor{newcolor}{rgb}{.8,.349,.1}
\usepackage[utf8]{inputenx}
\usepackage{graphicx}
\usepackage{amsmath,amssymb,mathrsfs}
\usepackage{fancyhdr}
\usepackage{calc}
\usepackage[T1]{fontenc}
\usepackage{mathptmx}
\usepackage{physics}
\usepackage{lineno}
\usepackage{tikz}
\usepackage{pgfplots}

\def\cp{$\cal{C}^+$}
\def\cm{$\cal{C}^-$}
\def\ww{\kappa}

\newcommand{\bes}{\begin{equation*}}
\newcommand{\ees}{\end{equation*}}
\newcommand{\beq}{\begin{equation}}
\newcommand{\eeq}{\end{equation}}
\newcommand{\bea}{\begin{eqnarray}}
\newcommand{\eea}{\end{eqnarray}}
\newcommand{\beas}{\begin{eqnarray*}}
\newcommand{\eeas}{\end{eqnarray*}}
\newcommand{\ds}{\displaystyle}
\newcommand{\gam}{\gamma}
\newcommand{\gamo}{\gamma_1}
\newcommand{\colb}[1]{\textcolor{black}{#1}}
\newcommand{\colt}[1]{\textcolor{black}{#1}}
\def\cp{$\cal{C}^+$}
\def\cm{$\cal{C}^-$}
\def\st{$\cal{S}$}

\begin{document}

\begin{frontmatter}

\title{Ordinary differential equations for the adjoint Euler equations}

\address[rvt1]{DAAA, ONERA, Universit\'e Paris Saclay, F-92322 Ch\^atillon, France}
\address[rvt2]{Centre Inria Universit\'e C\^ote d'Azur, Inria, 2004 Route des Lucioles, BP 93, 06902 Sophia Antipolis, France}

\author[rvt1]{Jacques Peter\corref{cor1}}
\ead{jacques.peter@onera.fr}
\cortext[cor1]{Corresponding author. Tel.: +33 1 46 73 41 84.}

\author[rvt2]{Jean-Antoine D\'esid\'eri}
\ead{jean-antoine.desideri@inria.fr}

\begin{abstract}
Ordinary Differential Equations are derived for the adjoint Euler equations
 firstly using the method of characteristics in 2D. For this system of 
 partial-differential equations, the characteristic 
 curves appear to be the streamtraces and the well-known \cp and \cm 
 curves of the 
  theory applied to the flow. The differential equations 
  satisfied along the streamtraces in 2D
   are then extended and demonstrated in 3D by linear combinations of the original 
   adjoint equations. 
These findings extend their well-known
 counterparts for the direct system, and should serve analytical and possibly numerical studies
of the perfect-flow model with respect to adjoint fields or sensitivity questions. 
Beside the analytical theory, the results are demonstrated 
by the numerical integration of the compatibility relationships
for discrete 2D flow-fields and dual-consistent adjoint fields over a very fine grid
about an airfoil.
\end{abstract}   


\end{frontmatter}

\section{\label{sec:l1_1} Introduction}
In 1988, Jameson derived the continuous adjoint equations associated with
the 2D and 3D Euler equations using general curvilinear coordinates \cite{Jam_88}.
With this landmark article, the fluid dynamics and aeronautical communities became better
aware of the potential of the adjoint approach for design, that is, the possibility to calculate
gradient information at a cost scaling with the number of functions to be differentiated, independently of the number of design parameters.
The equations in \cite{Jam_88} appeared to be a natural starting point for local optimizations  involving a large number of design variables by using adjoint gradients.
However, in that setting, the flow and the dual fields had to be calculated over a structured mesh. 
   
\colb{Nine years later}, Anderson and Venkatakrishnan \cite{AndVen_97,AndVen_99} and also Giles and
Pierce \cite{GilPie_97} derived the corresponding equations in Cartesian coordinates 
thus allowing the application of the continuous approach (sometimes referred to as the differentiate-then-discretize approach) on all types of meshes and, in particular, on unstructured meshes. 
%
For the sake of simplicity, we present here the two-dimensional case only in which 
the adjoint equations read
\beq
-A^T \frac{\partial \psi}{\partial x} - B^T \frac{\partial \psi}{\partial y}=0, \quad \mbox{in }\Omega
\mbox{ the fluid domain}
\label{e:adjcdom} 
\eeq
where $A$ and $B$ are the Jacobian matrices of the flux vectors $F_x$ and $F_y$ of the 
Euler equations in the $x$ and $y$ directions respectively:
\colb{
\beas
F_x =\left( \begin{array}{c}
   \rho u \\
   \rho u^2 + p  \\ 
   \rho u v  \\ 
   \rho u H  \\
\end{array}
\right)
\qquad
F_y =\left( \begin{array}{c}
   \rho v \\
   \rho u v   \\ 
   \rho v^2 + p  \\ 
   \rho v H  \\
\end{array}
\right), \nonumber 
\eeas
}
\colb{with $\rho$ the density, $(u,v)$ the velocity components, $p$ the static pressure and $H$ the total enthalpy}. In the most common case where the
quantity of interest \colb{(QoI)} is a line integral along the solid wall $\Gamma_w$, it can be shown easily that
the adjoint wall boundary condition is well-posed provided that the
function of interest depends only on the static pressure \cite{AndVen_99}.
In the classical case where the functional output of interest
is the force on $\Gamma_w$ projected in direction $\overline{d}$,  $J= \ds \int _{\Gamma_w} p (\overline{n} \cdot \overline{d})  ds$,
the wall boundary condition reads
\beq
   \overline{n} \cdot \overline{d} + \psi_2 n_x + \psi_3 n_y = 0 ~~~~ \textrm{on}~~ \Gamma_w.
 \label{e:adjcwal}
 \eeq
For the farfield of an external flow, as well as for the inlet and outlet of an internal flow,
the boundary conditions are derived from the theory of  local
one-dimensional characteristic decomposition \cite{AndVen_97, GilPie_01}. 
Here, the continuous adjoint Euler equations and
the associated boundary conditions are abbreviated as (AE).
Along with the growing use of the adjoint method for shape optimization, 
goal oriented mesh adaptation and also
meta-modelling, stability or control, great effort is being devoted 
to gain understanding in the mathematical properties
of the (AE) solutions. The main results are summarized here
before discussing the characteristic relations for the (AE) system.

After the derivation of the (AE) equations,
the first demonstrated property was also due to Giles and Pierce \cite{GilPie_97}: in the common
case where the function of interest is an integral along the wall, the authors proved that the
first and last components of
the adjoint vector $\psi$, associated with mass and energy conservation, 
satisfy $\psi_1 =  H \psi_4$.

Besides, the integration by parts yielding  (\ref{e:adjcdom}) is not valid in the entire fluid domain 
in the presence of flow discontinuities. After a series of works dealing with the quasi-1D Euler equations --
see \cite{GilPie_01,BaeCasPal_09} and references therein -- Baeza et al. presented the equations 
complementing (\ref{e:adjcdom}) along a shock line \cite{BaeCasPal_09} (denoted here $\Sigma$ as in the original reference).
The new equations are derived by introducing a complementary set of Lagrange multipliers,
multiplying the Rankine-Hugoniot conditions, viewed as constraints on $\Sigma$. 
Finally, the continuity of the adjoint field $\psi$ along  $\Sigma$ is established,
although $\nabla \psi$ may be discontinuous across $\Sigma$,
as well as $\psi$ over $\Gamma_w \cap \Sigma $, and a so-called internal boundary condition
is derived:
\bea
(\partial {\bf \psi}/\partial t).(F_x t_x + F_y t_y) = 0 \qquad \textrm{on} \qquad \Sigma    \label{e:adjcsi}   \\
\textrm{with}~ t ~\textrm{the unit vector tangent to }\Sigma  \nonumber 
\eea  
Coquel et al., Lozano and Renac \cite{CoqMarRai_18,Loz_18,JPRenLab_22} have derived additional relationships by using (\ref{e:adjcsi}),
 the jump operator applied to (\ref{e:adjcdom}) across $\Sigma$ and the Rankine-Hugoniot equations.

The fact that $\psi_1 =  H \psi_4$ can be proven simply by forming the linear
combination of the first three lines of system (\ref{e:adjcdom}) with coefficients $(1,u/2,v/2)$.
This yields
$\overline{U}.\nabla \psi_1- H \overline{U}.\nabla \psi_4 = 0  $ (with $\overline{U}=(u,v)$ 
the velocity vector).
Note that this was also derived in \cite{GilPie_97} by an approach based on physical source terms,
constituting an important analysis technique for the adjoint field of usual \colb{QoIs}. 
In particular, this method proved to be very fruitful to identify
the zones where numerical divergence of the adjoint vector is observed 
and mathematical divergence of
the solutions of (AE) is suspected. 
For the sake of clarity and brevity, we restrict the present discussion
to 2D flows about lifting airfoils, and to two of these zones, namely the stagnation streamline and the wall,
and to the lift and drag as functions of interest.
 
More precisely, Giles and Pierce \cite{GilPie_97} introduced four 
physical punctual source terms (or Green's functions in the classical mathematical vocabulary)
denoted here $\delta R^1$, $\delta R^2$, $\delta R^3$, $\delta R^4$.
These terms are added to the right hand-side of the linearised Euler equations 
and correspond respectively
to {\it (i)} a mass source at fixed stagnation pressure $p_0$ and enthalpy $H$ ; 
{\it (ii)} a normal force ; {\it (iii)} a change in $H$ at fixed $p$ and $p_0$; 
and {\it (iv)} a change in $p_0$ at fixed $p$ and $H$. 
They are linearly independent.
(We refer to the original reference for the detailed expression of these source terms.)
  The resulting changes in the \colb{QoI} $J$, $\delta J^l$,
  can be expressed as the integral over the domain of $\psi \delta R^l$
   that is, the value 
   at the source location
   since $\delta R^l$ is a Green's function. These source terms also admit a physical interpretation
   and their influence on the flow  can be understood in terms
of mechanical principles, and sometimes even quantified \colb{finally providing insight in the adjoint field \cite{GilPie_97}.}

  It has been observed 
  that the lift adjoint exhibits numerical divergence at the stagnation
  streamline and at the wall at subcritical flow conditions. Also the drag and lift adjoint of a transonic
  airfoil exhibit numerical divergence at the same locations if the foot of at least one shock wave
  is located strictly upwind the trailing edge -- see \cite{Loz_19,JP_20} and references therein.
  Reference \cite{JPRenLab_22} includes a careful verification of this physical perturbations approach
  applied to the discrete adjoint with a preliminary
   assessment of the consistency between the linear (discrete adjoint) and the non-linear (flow perturbation) 
 evaluations  of the $\delta J^l$.
  After this verification step, the non-linear perturbed flow approach has been used 
  (considering the physical source terms point of view prior to the classical adjoint) 
  and it appeared that:
  (a) $\delta R^4$ is the only source term causing a numerical divergence of $\delta J$ 
in the vicinity of the wall and stagnation streamline
; (b) in transonic condition, the numerical divergence of
   $\delta CLp^4$ and $\delta CDp^4$ in these zones is mainly due to the displacement of the shockfoot (or
   shock-feet if two shocks are not based at the trailing edge) ; (c) this numerical divergence is transferred
   to the adjoint components via the inverse matrix of the source terms ; (d) this does not necessarily
    prevent the numerical satisfaction
   of the adjoint lift- (resp. drag-) boundary condition at the wall (\ref{e:adjcwal}) as the
   calculation of $\psi_2 n_x + \psi_3 n_y$
   in this approach involves the product of $\delta CLp^4$ (resp. $\delta CDp^4$) by $(u n_x+v n_y)$.\\
   
  The method of characteristics for 2D inviscid supersonic flow is a classical method for deriving
  ordinary differential equations and, potentially, explicit algebraic relations satisfied along two
  families of curves, denoted \cp (left running with respect to (w.r.t.) a streamline) and \cm (right running).
  Here, we recall  the derivation of the continuous equations and study
  their counterparts for the (AE) equations.
   \\
 \colb{For the flow,} the method starting point is the Cauchy problem posed
   for the Euler equations: knowing the state variables along a fixed curve (L),
  is it possible to calculate their partial derivatives, in both space directions, for all points of (L)
  (this being a necessary condition for the flow calculation in the fluid domain) ?
  The flow variables in two neighboring points of (L), denoted here $a$ and $b$, 
  are linked by fluid dynamics equations and basic first order Taylor formulas.
  \colb{The more general presentations deal with rotational flows and, geometrically,
   both planar and axisymmetrical flows \cite{Fer_46,Sha_54,BonLun_89,Del_08}.
   The authors derive a system of
   equations for the derivatives of selected primitive variables}. 
  This system is linear in the unknown partial derivatives with non linear functions of the state
   variables as coefficients.
When its determinant is equal to zero, it cannot be solved. That is the case if $b$ is on the streamline of
  $a$ or\colb{, in case of a supersonic flow,} if the angle of $\overrightarrow{ab}=(dx,dy)$
   w.r.t. the streamline passing through
   $a$ is $\pm \sin^{-1}(1/M) $ ($M$ being the local Mach number).
   The specific curves (L)  where these conditions are satisfied for every points are, \colb{for all Mach numbers}, the streamtraces, \colb{and,
   where the flow is supersonic,} the
   so-called \cp (left running curves w.r.t. the streamtraces with angle $\sin^{-1}(1/M) $) and the \cm
   (right running curves the  w.r.t. streamtraces  with angle $-\sin^{-1}(1/M) $).
   Along these curves, the physical existence and boundedness of the vector of unknowns
   allows to conclude from the nullity of the determinant in the denominator of Cramer's formulas, to the
   nullity of the determinants appearing in the numerators, and this, for all the
   variables. 
 The classical computational method for supersonic flow \cite{Fer_46,Sha_54,LipRos_56,BonLun_89,And_03,Del_08} is supported by the
corresponding differential forms valid along the \cp and \cm curves, and \colb{the property of constant
total enthalpy and constant entropy} along the streamtraces.
 If the flow is irrotational, homoenthalpic and homoentropic,
 simpler equations are derived for the velocity magnitude and the velocity 
 angle \cite{Sha_54,Sea_54,LipRos_56} or the  velocity potential \cite{And_03} and 
 the differential forms satisfied along the \cp and \cm curves may be integrated.
 This permits to establish the well-known equations
\beq   
 k^- = \phi + \nu(M)~~ \textrm {is constant along a}~~ {\cal C}^- ~~~~~~~~~~~~~~~~~
      k^+ = \phi - \nu(M)~~ \textrm {is constant along a}~~ {\cal C}^+ 
\label{e:kpkm}
\eeq
in which $\phi$ is the streamline angle, $ \nu(M)$ the Prandtl-Meyer function,
$$ 
\phi=\tan^{-1} (v/u) ,\quad
\nu(M) =\sqrt{\frac{\gamma+1}{\gamma-1}}
      \tan^{-1} (\sqrt{\frac{\gamma-1}{\gamma+1} (M^2 -1)}) -\tan^{-1}(\sqrt{M^2-1}),
 $$ 
  and $\gamma$ the ratio of specific heats ($\gamma=\frac{7}{5}$ for diatomic perfect gas).
  \colb{Finally, let us recall that the integration of the
  corresponding ODEs along the trajectories results in the 
   property of constant enthalpy and constant entropy.}
 \\ 
 Besides, Bonnet and Luneau indicate that the mechanical equations posed in $a$ may be expressed in the
 Cartesian frame of reference rather than in the usual local frame derived from the local velocity
 \cite{BonLun_89}. Note that the 2D characteristic equations could also be calculated in the
 Cartesian frame and in an inexpert way, without
   taking advantage of the known properties of the streamtraces.
   Then the following the $8\times8$ linear system relating the derivatives of
   the conservative variables  would be solved:
  \begin{equation} 
\begin{bmatrix}
dx & 0  & 0  & 0  & dy & 0  & 0  & 0 \\
0  & dx & 0  & 0  & 0  & dy & 0  & 0 \\
0  & 0  & dx & 0  & 0  & 0  & dy & 0 \\
0  & 0  & 0  & dx & 0  & 0  & 0  & dy \\
   &    &    &    &    &    &    &   \\
   &    &    &    &    &    &    &   \\
   &    & A  &    &    &    & B  &   \\
   &    &    &    &    &    &    &   
\end{bmatrix}
\begin{bmatrix}
(\partial \rho /\partial x) \\
(\partial \rho u/\partial x) \\
(\partial \rho v/\partial x) \\
(\partial \rho E/\partial x) \\
(\partial \rho /\partial y) \\
(\partial \rho u/\partial y) \\
(\partial \rho v/\partial y) \\
(\partial \rho E/\partial y) 
\end{bmatrix}
=
\begin{bmatrix}
  \rho^b - \rho^a\\
  \rho u^b - \rho u^a\\
  \rho v^b - \rho v^a\\
  \rho E^b - \rho E^a\\
 0. \\
 0. \\
 0.\\
 0.
\end{bmatrix}
\label{e_direct}
  \end{equation}
The starting point of our analytical development resides in the observation that 
(\ref{e_direct}) and the corresponding linear system
for the (AE) equations, (\ref{e_cbase}), have the same determinant. 
From this observation, the adjoint Euler characteristics equations
are established in Sec. 2.
The theoretical findings are \colb{linked with former researches and} illustrated by numerical computational solutions 
over a very fine grid in Sec. 3.
Conclusions are drawn in Sec. 4.
%
\section{\label{sec:l1_2} Adjoint characteristic equations for 2D supersonic flow}
%
The method exposed in \cite{And_03} (resp. \cite{Del_08}) for potential (resp. general) inviscid
flow has served as a guideline to our derivation for the adjoint system. For all our calculations,
 we assume an ideal gas law for
 the static pressure $p = (\gamma - 1)\rho e = (\gamma - 1)(\rho E - 0.5  \rho || \overline{U} || ^2)$ with a
  constant $\gamma$.
%
\subsection{Problem statement}
%
Given two fixed close points in the supersonic zone, $a$ and $b$, 
is it possible to estimate 
$(\partial \psi / \partial x), (\partial \psi / \partial y)$ from the local value of the flow field and
($\psi_a$,$\psi_b$) ? This question is the starting point of the method of characteristics in which
specific lines are identified along which this problem is ill-posed, and particular  
ordinary differential equations are satisfied.
Let us denote 
$\overrightarrow{ab}=(dx,dy)$  and first  assume that $dx \neq 0$. 
By definition of differential forms, and in view of the adjoint system
(\ref{e:adjcdom}), the following holds 
\begin{equation} 
\begin{bmatrix}
dx & 0  & 0  & 0  & dy & 0  & 0  & 0 \\
0  & dx & 0  & 0  & 0  & dy & 0  & 0 \\
0  & 0  & dx & 0  & 0  & 0  & dy & 0 \\
0  & 0  & 0  & dx & 0  & 0  & 0  & dy \\
   &    &    &    &    &    &    &   \\
   &    &    &    &    &    &    &   \\
   & -  & A^T  &    &    & -  & B^T  &   \\
   &    &    &    &    &    &    &   
\end{bmatrix}
\begin{bmatrix}
(\partial \psi_1/\partial x) \\
(\partial \psi_2/\partial x) \\
(\partial \psi_3/\partial x) \\
(\partial \psi_4/\partial x) \\
(\partial \psi_1/\partial y) \\
(\partial \psi_2/\partial y) \\
(\partial \psi_3/\partial y) \\
(\partial \psi_4/\partial y) 
\end{bmatrix}
=
\begin{bmatrix}
 d\psi_1 \\
 d\psi_2 \\
 d\psi_3 \\
 d\psi_4 \\
 0 \\
 0 \\
 0\\
 0
\end{bmatrix}
\label{e_cbase}
\end{equation}
in which by neglecting second-order terms in space:
$(d\psi_1, d\psi_2, d\psi_3, d\psi_4)=(\psi^b_1-\psi^a_1, \psi^b_2-\psi^a_2, \psi^b_3-\psi^a_3, \psi^b_4-\psi^a_4)$.
The determinant of the linear system is evidently
\begin{equation*} 
\begin{vmatrix}
      &    &    &    &    &    &   \\
      & dx I&    &    &  & dy I  &   \\
      &    &    &    &    &    &   \\
      &    &    &    &    &    &   \\
      &    &    &    &    &    &   \\
      & -A^T  &    &    &  & -B^T  &   \\
      &    &    &    &    &    &   
\end{vmatrix}
= 
\begin{vmatrix}
   &    &    &    &    &    &   \\
   & dx I &  &    & 0  &   &   \\
   &    &    &    &    &    &   \\
   &    &    &    &    &    &   \\
   &    &    &    &    &    &   \\
   & -A^T&   &    & -B^T+dy/dx A^T  &  & \\
   &    &    &    &    &    &   
\end{vmatrix}
= dx^4\vert - B^T + dy/dx A^T \vert  = \vert -dx B^T + dy A^T \vert 
\end{equation*}
Of course,  $\vert -dx B^T + dy A^T \vert$ is equal to  $\vert -dx B + dy A \vert$ and the value of this determinant is known
from the eigenvalues of the matrix:
$$D = \vert -dx B + dy A \vert = (-v ~dx  + u ~ dy)^2 (-v ~dx  + u ~ dy  + c ~ds)(-v ~dx  + u ~ dy  - c ~ds),   $$
in which
$$  c= \sqrt{\frac{\gamma p}{\rho}},\quad ds = \sqrt{dx^2 + dy^2}   .$$
Similarly to the flow derivatives reconstruction \cite{And_03,Del_08}, the problem of adjoint derivatives
 reconstruction in a supersonic zone is ill-posed along the same three families of curves
\bea
-v ~dx  + u ~ dy  &=& 0 \qquad \qquad  \textrm{\colb{\st~ streamtraces \qquad (all Mach numbers)}}   \label{e_ctraj}\\
-v ~dx  + u ~ dy  + c ~ds &=& 0 \qquad \qquad \textrm{\cm characteristics \quad\colb{(supersonic flow only)}}\label{e_ccp}\\
-v ~dx  + u ~ dy  - c ~ds &=& 0 \qquad \qquad \textrm{\cp characteristics \quad\colb{(supersonic flow only)}}\label{e_ccm}
\eea
\begin {figure}[htbp]
  \begin{center}
	  \includegraphics[width=0.9\linewidth]{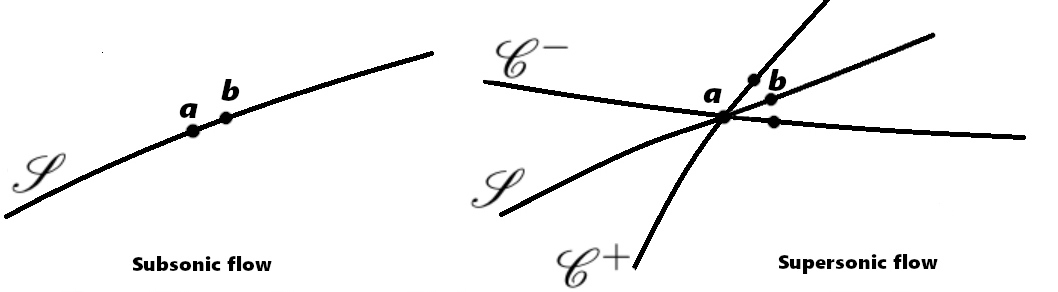}
	  \caption{\colb{\st, \cp and \cm curves according to local Mach number. Neighboring $a$ and $b$ points along these curves}}
  \label{f:stcpcmcurves}
  \end{center}
\end {figure} 
Classically, the method of characteristics uses the ill-posedness of (\ref{e_cbase}) in the following way: 
along the curves  defined by equations (\ref{e_ctraj}), (\ref{e_ccp})
or (\ref{e_ccm}), not only the denominator appearing in the Cramer formulas 
applied to the linear equations (\ref{e_cbase}) is equal to zero, but the numerators giving
the eight components of   $(\partial \psi / \partial x)$ and $ (\partial \psi / \partial y)$ must also be equal
to zero for the fractions not to be singular. This (somehow paradoxical) technique
allows the derivation of equations (\ref{e:kpkm}). It is derived here for the (AE) system by
 analysing the set of linear equations (\ref{e_cbase}).
%
%
\subsection{\colb{Null differential forms in the adjoint variations along trajectories and characteristics}}
%
The transposed of the Euler flux Jacobian matrices in $x$ and $y$ direction read 
$$
A^T = 
\begin{bmatrix} 
	0  &  \gamma_1 E_c -u^2  & -uv & (\gamma_1 E_c - H)u\\
	1  & (3-\gamma)u &  v  & H-\gamma_1 u^2\\
	0  & -\gamma_1 v &  u  & -\gamma_1 uv \\
	0  & \gamma_1 &  0  & \gamma u\\
\end{bmatrix}
~~~~~~~~
B^T = 
\begin{bmatrix} 
	0  & -uv & \gamma_1 E_c -v^2 & (\gamma_1 E_c - H)v\\
	0  &  v  & -\gamma_1 u & -\gamma_1 u v \\
	1  &  u  & (3-\gamma)v  & H -\gamma_1 v^2\\
	0  &  0  & \gamma_1       & \gamma v \\
\end{bmatrix}
$$
in the usual notations in aerodynamics and $\gamma_1 = \gamma-1$. 
Let
$$ t = \frac{dy}{dx} ,\quad \ww = u t -v,$$
and also introduce the following notations for the column vectors of the transposed Jacobian matrices:
 $ A^T = [A_1|A_2|A_3|A_4] $, $ B^T = [B_1|B_2|B_3|B_4] $. 
Before presenting the results, the principle of the calculation is recalled
in one of the cases that leads to the simplest calculations: 
the definition of $(\partial \psi_4/\partial x)$ along the curves
where  (\ref{e_ctraj}), (\ref{e_ccp}) or (\ref{e_ccm}) is satisfied
(that is,  the streamtraces, the \cm or \cp characteristic) 
requires that, along these curves

\begin{equation*} 
\begin{vmatrix} 
  ~dx  & 0  & 0  &  d\psi_1   & dy   &   0  &  0 & 0~ \\
  ~0   & dx & 0  &  d\psi_2  &  0   & dy   &  0 & 0 ~\\
  ~0   & 0  & dx &  d\psi_3  &  0   &   0  & dy & 0 ~\\
  ~0   & 0  & 0  &  d\psi_4  &  0   &   0  & 0  & dy ~\\
  ~|   & |   &  |  &   |      &  |    &  |    &  |  & |~\\
  ~-A_1 & -A_2 &  -A_3  & 0      &  -B_1  & -B_2 & -B_3 & -B_4 ~\\ 
  ~|   & |   &  |  &   |      &  |    &  |    & |   & |~\\
\end{vmatrix}
 = 0
\end{equation*}
The determinant is expanded along the fourth column and the following notations are used
\beq
 -C_{4x}^1 d\psi_1 + C_{4x}^2 d\psi_2 - C_{4x}^3 d\psi_3 + C_{4x}^4 d\psi_4 = 0      
\label{e_c4x}
\eeq
in which, for example
\begin{equation*}
  C_{4x}^1=
\begin{vmatrix} 
  ~0   & dx  & 0    &  0   &   dy   &  0  & 0~ \\
  ~0   & 0   & dx   &  0   &   0    & dy  & 0 ~\\
  ~0   & 0   & 0    &  0   &   0    &  0  & dy ~\\
  ~|   & |   &  |   &  |    &  |    &  |  & |~\\
  ~-A_1 & -A_2 &  -A_3 &  -B_1  & -B_2 & -B_3 & -B_4 ~\\ 
  ~|   & |   &  |   &  |    &  |    & |   & |~\\
\end{vmatrix}
=
\begin{vmatrix} 
  ~0   &  dx  & 0   &  0   &   0  &  0 & 0~ \\
  ~0   & 0    & dx   &  0   &   0  &  0 & 0 ~\\
  ~0   & 0    & 0  &  0   &   0  &  0 & dy ~\\
  ~|   & |   &  |   &  |    &  |     &  |  & |~\\
  ~-A_1 & -A_2 &  -A_3 &  -B_1   & -B_2+t A_2 & -B_3 +t A_3 & -B_4 ~\\ 
  ~|   & |   &  |   &  |    &  |    & |   & |~\\
\end{vmatrix} 
\end{equation*}
Finally
$$  C_{4x}^1=  dx^2 dy ~ \vert~~ -A_1~~ -B_1 ~~ (-B_2+ t A_2)~~ (-B_3+t A_3) ~~\vert $$ 
The determinant of this $4 \times 4$ matrix is easily calculated thanks to the simplicity of the first two columns $A_1$ and $B_1$. The
final result is
$$  C_{4x}^1=   dx^2 dy ~ \gamma_1 ~\ww ~(u+v t)  $$
We emphasise  that at this stage no assumption is made on the value of $t=dy/dx$ w.r.t. 
the velocity vector $(u,v)$. In particular $t$ is not assumed to be the
tangent of the  angle of the velocity w.r.t. the $x$ axis and $\ww$ is not assumed to be zero. 
This is mandatory
to derive relations that can be used for all three types of specific curves 
and also to account for the multiplicity
of the eigenvalue $ (-v ~dx  + u ~ dy) $ along streamtraces.
The other terms of the differential form of interest read
\beas
 C_{4x}^2 &=&    - dx^2 dy ~ \vert~~ -A_2~~ (-B_1+t A_1) ~~ -B_2~~ (-B_3+t A_3) ~~\vert =  - dx^2 dy ~ \gamo ~ \ww ~(u^2 + v^2) \\
 C_{4x}^3 &=&      dx^2 dy ~ \vert~~ -A_3~~ (-B_1+ t A_1) ~~ (-B_2+ t A_2)~~ -B_3 ~~\vert =  dx^2 dy ~ \gamo ~\ww~ t~ (u^2 + v^2)\\
 C_{4x}^4 &=&      dx^3    ~ \vert~~ (-B_1+t A_1) ~~ (-B_2+ t A_2)~~ (-B_3+t A_3) ~~-B_4~~\vert \\
         &=&    - dx^3 \ww  \left( (\gamo+\gam t^2) u^2 v-2 u v^2 t + \gamo H \ww  + (\gam+\gamo t^2) v^3 -\gamo v~ (1+t^2)~ E_c \right)
 \eeas
\colb{The explicit expression of (\ref{e_c4x}), the necessary condition for the boundedness of 
$(\partial \psi_4/\partial x)$, hence reads}
\bea
&-& dx^2 dy~\gamma_1 ~\ww ~(u+v t)~d\psi_1 -  dx^2 dy~\gamo ~ \ww ~(u^2 + v^2) d\psi_2 -dx^2 dy~\gamo ~\ww~ t~ (u^2 + v^2) d\psi_3  \nonumber  \\ 
 &-& dx^3 \ww  \left( (\gamo+\gam t^2) u^2 v-2 u v^2 t + \gamo H \ww  + (\gam+\gamo t^2) v^3 -\gamo v~ E_c~ (1+t^2)\right) d\psi_4 = 0  \label{e_c4xa} 
\eea
\colb{Assuming that $dx\neq0$, this equation may be
further simplified for the \cm and \cp for which $\ww\neq0$:}
\bea
&& \gamma_1  ~t (u+v t)~d\psi_1 + \gamo  ~t (u^2 + v^2) ~ d\psi_2 +~\gamo ~ t^2~ (u^2 + v^2) ~ d\psi_3  \nonumber  \\ 
 &+&   \left( (\gamo+\gam t^2) u^2 v-2 u v^2 t + \gamo H \ww  + (\gam+\gamo t^2) v^3 -\gamo v~ E_c~ (1+t^2)\right) d\psi_4 = 0.   \label{e_c4xb} 
\eea
\colb{As $(-v ~dx  + u ~ dy)=\ww~dx$ has a multiplicity of two in the determinant of (\ref{e_cbase}), equation (\ref{e_c4xb}) is also
needed for the existence of $(\partial \psi_4/\partial x)$ and hence true for neighboring points $a$ and $b$ along the same \st curves. (This point is detailed in \S 2.3.)\\
For the sake of clarity, the results of the corresponding calculations for the existence of the seven other partial derivatives along the \st, \cm and \cp curves, 
and the properties of the $C^i_{jx}$  $C^i_{jy}$ coefficients are presented in Appendix A and B. Only the counterparts of equation (\ref{e_c4xb}), for the existence of $(\partial \psi_1/\partial x)$, $(\partial \psi_2/\partial x)$, $(\partial \psi_3/\partial x)$,
 $(\partial \psi_1/\partial y)$, $(\partial \psi_2/\partial y)$, $(\partial \psi_3/\partial y)$ and
  $(\partial \psi_4/\partial y)$ are presented hereafter in this order:}
\bea
& &  ( (2 t u + (t^2-1) v)(\gamo H + \gamo E_c + \gam v \ww) -(u+vt)(\gamo t H + \gamo u \ww + \gam v \ww t) d\psi_1   \nonumber  \\
&+&  ~ \gamo ~  t~ (u^2 + v^2) ~H ( d\psi_2 + t d\psi_3) +\gamo ~  t~ (u + v t) ~H^2  d \psi_4  = 0  \label{e_c1xb} 
\eea
\bea
& &  t~(\gamo H + \gamo E_c + \gam v \ww) ~(d\psi_1+H d\psi_4)    \nonumber  \\
&-& ( \gamo u^2 v -u v^2 t + \gam v^3  - \gamo (u t + v) ( E_c + H))~ (d\psi_2 + t d\psi_3)  = 0 \label{e_c2xb}  
\eea
\bea
&&  t ~(t \gamo H + t \gamo E_c - \gam u \ww ) (d\psi_1 + H d\psi_4)
+ t~( (\gam+1) u^2 v - t \gamo u v^2 + (v + u t) ( \gamo E_c +\gamo H -\gam u^2))d\psi_2  \nonumber  \\
&+& ~ ( (t v^2 - u \ww)(-\gamo u +(\gam+1) t v) - (u t-(1+2 t^2)  v) (\gamo E_c + \gamo H -\gam v^2)) d\psi_3  = 0  \label{e_c3xb}  
\eea
\bea
 & &  \left( (\gamo+\gam t^2) u^3 -2 u^2 v t - \gamo t \ww H  + (\gam+\gamo t^2) u v^2 -\gamo u~ (1+t^2)~ E_c \right) d\psi_1   \nonumber  \\
&+&  ~ \gamo ~ (u^2 + v^2) ~H ( d\psi_2 + t d\psi_3)   
+   ~\gamo  ~  (u + v t) ~H^2  d \psi_4  = 0  \label{e_c1yb}
\eea
\bea
& &  (\gamo H + \gamo E_c + \gam v \ww) ~(d\psi_1+H d\psi_4) \nonumber \\
&+&  ( (v\ww+u^2)(vt-(\gamo+\gam t^2)u)-(vt-(2+t^2)u)(\gamo E_c + \gamo H +\gam v\ww) ) d\psi_2    \nonumber  \\
&-& ( \gamo u^2 v -u v^2 t + \gam v^3  -\gamo (u t + v) ( E_c + H)) d \psi_3  = 0 \label{e_c2yb}  
\eea
\bea
&&   ~(t \gamo H + t \gamo E_c - \gam u \ww ) (d\psi_1 + H d\psi_4) \nonumber \\
&+& ( (\gam+1) u^2 v - t \gamo u v^2 + (v + u t) ( \gamo E_c +\gamo H -\gam u^2))(d\psi_2 + t d\psi_3)  =0  \label{e_c3yb} 
\eea
\bea
&&  ~\gamma_1 ~(u+v t)~d\psi_1 +  \gamo  ~(u^2 + v^2) (d\psi_2 + t d\psi_3 ) \nonumber  \\ 
 &+&   \left( (\gamo+\gam t^2) u^3 -2 u^2 v t - \gamo t \ww H + (\gam+\gamo t^2) u v^2 -\gamo u~ E_c~ (1+t^2)\right) d\psi_4 = 0 \label{e_c4yb} 
\eea
\subsection{Ordinary differential equations for the adjoint along the streamtraces \colb{\st}}
%
The trajectories are one of the families of specific curves for the gradient calculation problem
 (\ref{e_cbase}). Along these curves
$ u~dy-v ~dx = 0 $ is a zero of the denominator of Cramer's formulas
with multiplicity two.
Let us first assume that point $a$ is fixed and point $b$ is very close to  ${\cal S}_a$, the streamtrace passing through
 $a$  but not on this curve. The first-order expression of $(\partial \psi_1 / \partial x)$
is
$$ \frac{\partial \psi_1}{ \partial x} = \frac{C^1_ {1x} d\psi_1 - C^2_ {1x} d\psi_2+ C^3_ {1x} d\psi_3 - C^4_ {1x} d\psi_4}
 { (-v ~dx  + u ~ dy)^2 (-v ~dx  + u ~ dy  + c ~ds)(-v ~dx  + u ~ dy  - c ~ds)}  $$
 Actually $\ww dx = u~dy-v~dx $  is a factor of all four coefficients $C^1_ {1x}$, $C^2_ {1x}$, $C^3_ {1x}$ and $ C^4_ {1x}$. We denote by
 $\bar{C}^l_ {mx}$ the coefficients obtained by removing the $ (\ww~dx)$ factor from the corresponding
 $C^l_ {mx}$. Obviously 
 $$ \frac{\partial \psi_1}{ \partial x} = \frac{\bar{C}^1_ {1x} d\psi_1 - \bar{C}^2_ {1x} d\psi_2+ \bar{C}^3_ {1x} d\psi_3 - \bar{C}^4_ {1x} d\psi_4}
 { (-v ~dx  + u ~ dy) (-v ~dx  + u ~ dy  + c ~ds)(-v ~dx  + u ~ dy  - c ~ds)}  $$
 If point $b$ is moved closer and closer to ${\cal S}_a$ $ (-v ~dx  + u ~ dy) \rightarrow 0 $,  so that the boundedness of   $(\partial \psi_1/ \partial x)$ requires that
$$ \bar{C}^1_ {1x} d\psi_1 - \bar{C}^2_ {1x} d\psi_2+ \bar{C}^3_ {1x} d\psi_3 - \bar{C}^4_ {1x} d\psi_4 = 0 ~~~~~~~~~~~~ \textrm{on} ~~~{\cal S}_a   $$
   This expression and its counterparts for the other derivatives    $(\partial \psi_2/ \partial x)$ ...  $(\partial \psi_4/ \partial y)$ have to be satisfied
        for all trajectories. How many of these eight differential forms are independent ? If $ \ww~ dx = 0 $, then
\beq
\bar{C}^1_{1x} = - t ~\bar{C}^1_{1y} ~~~~~ \bar{C}^2_{2x} = - t ~\bar{C}^2_{2y} ~~~~ \bar{C}^3_{3x} = - t ~\bar{C}^3_{3y} ~~~~  \bar{C}^4_{4x} = - t ~\bar{C}^4_{4y}
\label{e_prop5}
\eeq
as in this case 
$$ ~\bar{C}^1_{1x} = \bar{C}^4_{4x} = -0.5 ~t \gamo dx^2 (1+t^2)^2 u^3  ~~~~~~ ~\bar{C}^1_{1y} = \bar{C}^4_{4y} = 0.5 ~\gamo dx^2 (1+t^2)^2 u^3  $$
$$ ~\bar{C}^2_{2x} =  2~ dx^2 \gamo t u H ~~~~~ \bar{C}^2_{2y} = - 2~ dx^2 \gamo u H ~~~~~~~~~~~~~~~
   \bar{C}^3_{3x} =  2~ dx^2 \gamo t^3 u H ~~~~~ \bar{C}^3_{3y} = - 2~ dx^2 \gamo t^2 u H. $$
   Relations (\ref{e_cxy1}) to (\ref{e_cxy4}) are valid for the $\bar{C}$ coefficients
   (as they stand whatever the values of $w$ and $dx$, they may be simplified
   by $w d x$). In the specific case where   $\ww =0~$ they are completed by (\ref{e_prop5}).
    Equations \{(\ref{e_c1yb}),(\ref{e_c2yb}),(\ref{e_c3yb}),(\ref{e_c4yb})\}
     (necessary for the boundedness   
     of the $\partial \psi^l / \partial y$)
      and \{(\ref{e_c1xb}),(\ref{e_c2xb}),(\ref{e_c3xb}),(\ref{e_c4xb})\}
     (same for $\partial \psi^l / \partial x$) are then proportional by a $(-t)$ factor. 
   Considering the range of the set of the eight differential forms, it appears that 
   one of these two sets of four equations need be accounted for.
   
   The relations stemming from the existence of the $\partial x$ partial derivative
   are retained. Equation (\ref{e_c1xb}), required for the definition of $(\partial \psi_1/\partial x)$, is further simplified using the specific properties
    of a trajectory ($\ww=0~~ v=ut$):
 $$
   0.5 t \gamo (1+t^2)^2 u ^3 d\psi_1 +  ~ \gamo ~ t~ (1 + t^2) u^2 ~H ( d\psi_2 + t d\psi_3)  +   ~\gamo ~  t ~ (1 + t^2) u ~H^2  d \psi_4  = 0  $$
   $$  0.5  (1+t^2) u ^2 d\psi_1 +   u ~H ( d\psi_2 + t d\psi_3)  +   ~H^2  d \psi_4  = 0  $$
   For the streamstraces,
 the equation finally derived from the existence of $(\partial \psi_1/\partial x)$ is
\beq
   E_c ~d\psi_1 +   H ( u~ d\psi_2 + v~ d\psi_3)  +   ~H^2  d \psi_4  = 0 
   \label{e_atraj1}
\eeq
      Note that we have assumed that $u\neq0$ and $dx\neq0$ to perform the calculations but finally obtained an expression that is also well-defined in this
   specific case. Equation (\ref{e_c2xb}), is further simplified for the motion along a trajectory:
  $$   ~t~(\gamo H + \gamo E_c ) ~(d\psi_1+H d\psi_4) +2 \gamo~ u~ t~ H~ (d\psi_2 + t d\psi_3)  = 0 $$
  $$    (H + E_c ) ~(d\psi_1+H d\psi_4) +2  H~ (u~ d\psi_2 + v~ d\psi_3)  = 0 $$
   Using the first relation  and simplifying by $H$, we get
      $$    H (d\psi_1+H d\psi_4) +2  H~ (u~ d\psi_2 + v~ d\psi_3) +E_c H d\psi_4 -  H ( u d\psi_2 + v d\psi_3)  -  ~H^2  d \psi_4   = 0 $$
\beq
 \textrm{and finally} ~~~~~~~ d\psi_1 + u ~ d\psi_2 +  v ~ d\psi_3 +  E_c d\psi_4 = 0 
\label{e_atraj2}
\eeq
Simplifying equation (\ref{e_c3xb}) for trajectories yields  
$$  (H + E_c ) (d\psi_1 + H d\psi_4) +2 H ( u d\psi_2+  v d\psi_3) = 0, $$
that had already been derived.
Using  $\ww=0$ and $v=ut$, the fourth relation (first simplified by a $\ww dx$ factor) gives
$$ \gamo ~t u (1+ t^2)~d\psi_1 + \gamo  ~t u^2 (1 + t^2) ~ d\psi_2 +~\gamo  t^2~u^2 (1 + t^2) ~ d\psi_3 +  0.5 t \gamo (1+t^2)^2 u^3 d\psi_4 = 0 $$
$$  \textrm{or}~~~~~~~~~~~~~    d\psi_1 + u ~ d\psi_2 +  t u ~ d\psi_3 +  0.5 (1+t^2) u^2 d\psi_4 = 0 $$
$$  \textrm{and finally}~~~~~~~~~~~~~~~      d\psi_1 + u ~ d\psi_2 +  v ~ d\psi_3 +  E_c~ d\psi_4 = 0 $$
This equation is similar to (\ref{e_atraj2}).
These differential forms along the trajectories ${\cal S}$ may be turned in differential equations. The natural variable w.r.t. which differentiate, is the curvilinear abscissa
along the streamtraces, $s$, increasing in the direction opposite to the motion (as this is the direction of the adjoint transport of information from the support
of the function of interest). The final equations along the  ${\cal S}$ curves are then
\bea
  E_c ~\frac{d\psi_1}{ds} +   H ( u~ \frac{d\psi_2}{ds} + v~ \frac{d\psi_3}{ds})  +   ~H^2  \frac{d \psi_4}{ds}  = 0   \label{e_atraja} \\
  \frac{d\psi_1}{ds} + u ~ \frac{d\psi_2}{ds} +  v ~ \frac{d\psi_3}{ds} +  E_c \frac{d\psi_4}{ds} = 0    \label{e_atrajb}
  \eea
   Although the calculations of this section are not very complex,
        they have been checked with the computer algebra software Maple.
  The independent variables of the formal calculations are $(u,v)$ also $M$ the Mach number
  (that allows the calculation of the speed of sound $c$ and then the total enthalpy $H$) and $\gamma$.
  The set of four differential
  forms \{(\ref{e_c1yb}), (\ref{e_c2yb}), (\ref{e_c3yb}), (\ref{e_c4yb})\} was associated with
 \{(\ref{e_atraj1}),(\ref{e_atraj2})\} in a $6 \times4$ matrix that was again found
  to have rank two. 
  %
\subsection{Ordinary differential equation for the adjoint along the \cp and \cm characteristics }
%
Let us first note that the determinant in the denominator of the Cramer formulas for (\ref{e_cbase})
may be calculated as
$$ D = dx ~C^{1}_{1x} + dy ~C^{1}_{1y} $$
developing $D$ along the first line.
Doing the same along the second, third and fourth lines yields
$$  D = dx~ C^{2}_{2x} + dy~ C^{2}_{2y} = dx~ C^{3}_{3x} + dy~ C^{3}_{3y} = dx~ C^{4}_{4x} + dy~ C^{4}_{4y}   $$
As $D=0$ along the  \cp  and \cm characteristics, equations  (\ref{e_cxy1}) to (\ref{e_cxy4}) are completed by
$$ C^1_{1x} = - t ~C^1_{1y} \qquad C^2_{2x} = - t ~C^2_{2y} \qquad C^3_{3x} = - t ~C^3_{3y} \qquad  C^4_{4x} = - t ~C^4_{4y} $$
so that the differential forms stemming from the boundedness of $(\partial \psi^l/\partial x)$ and their
counterparts for $(\partial \psi^l/\partial x)$ are proportional. (Incidentally, note that this argument
 may have been used also in the previous subsection.)
\\
Numerical tests indicate that the four differential forms associated with the existence of (choosing
the second set) $(\partial \psi_1/\partial y)$, $(\partial \psi_2/\partial y)$, $(\partial \psi_3/\partial y)$,
$(\partial \psi_4/\partial y)$ are proportional but the corresponding calculations are much more 
complex than in the previous subsection as the expression of $t$ is now:
\beq
\colt{t^\pm} = \tan(\phi + \beta)~~~~ \textrm{with}~~~ \tan \phi = \frac{v}{u} ~~~ \sin \beta = \pm \frac{1}{M}
\label{etapb}
\eeq
\colt{with, of course $ \sin \beta = 1/M $ for the \cp and
$ \sin \beta = -1/M$ for the \cm curves.}\\
 With this expression of $t$ and $\ww$ not being equal to zero, searching the rank of \{(\ref{e_c1yb}),
  (\ref{e_c2yb}), (\ref{e_c3yb}), (\ref{e_c4yb})\} is much more difficult.
  However, the task was successfully accomplished with the assistance of Maple,
  using once again $(u,v,M,\gamma)$ as independent formal variables. 
  Rank one was indicated by Maple and the correspondence with the
  counterpart characteristics for the flow seemed very sound. 
  Knowing this result from formal calculation, its demonstration by hand
  was searched for.
\\
First, the value of $t$ is calculated along the \cp and \cm curves.
For $\beta$ in $[-\pi/2,\pi/2]$ (\ref{etapb}) yields
\beq
t^\pm =  \frac{\ds\frac{v}{u}\pm\ds\frac{1}{\sqrt{M^2-1}}}{1\mp\ds\frac{v}{u}\ds\frac{1}{\sqrt{M^2-1}}} = \frac{uv M^2 \pm (u^2 + v^2)\sqrt{M^2 -1} }{u^2 (M^2 -1)-v^2}
= \frac{uv \pm c^2\sqrt{M^2 -1} }{u^2-c^2}= \frac{uv \pm c\sqrt{u^2 + v^2 -c^2} }{u^2-c^2}
\label{eteps}
\eeq
We then note that if the differential forms  (\ref{e_c1yb}) and (\ref{e_c4yb}) were proportional, from
 the expressions of $C^2_{1y}, C^3_{1y}, C^2_{4y} $ and $ C^3_{4y}$, the ratio
 between their terms would be $(-H)$. It is then easily checked that the two non-trivial conditions
 for proportionality, 
 $C^1_{1y} = -H C^4_{1y}$ and $ C^1_{4y} = -H C^4_{4y}$, are equivalent to a single equality
\beq
\gamma_1 ~(u+v t) H =  \left( (\gamo+\gam t^2) u^3 -2 u^2 v t - \gamo t \ww H  + (\gam+\gamo t^2) u v^2 -\gamo u~ (1+t^2)~ E_c \right)
\label{e_cond1eq4}
\eeq
 Wherever $u\neq0$, this condition is equivalent to
 \beq
  \gamma_1 ~(1+t^2)~ H =  \gamo~(1+ t^2)~E_c + (tu -v)^2
  \label{e_deg2}
 \eeq
 that is, precisely, the degree two equation which roots are the values of $t$ along the \cp and \cm curves (\ref{eteps}).
 Along these curves, using (\ref{e_cond1eq4}) that is now an established property along the \cp and \cm, these two differential forms
  may be simplified as
  \beq
  (u+v t^\pm)~d\psi_1 +  (u^2 + v^2) (d\psi_2 + t^\pm d\psi_3 ) + H  (u+v t^\pm)~d\psi_4 = 0,
   \label{e_c14y} 
   \eeq
or, under the form of an ordinary differential equation,
    \beq
  (u+v t^\pm)~\frac{d\psi_1}{ds} +  (u^2 + v^2) (\frac{d\psi_2}{ds} + t^\pm \frac{d\psi_3}{ds} ) + H  (u+v t^\pm)~\frac{d\psi_4}{ds} = 0.
   \label{edo_c14y} 
   \eeq
  Comparing equations (\ref{e_c3yb}) -- expressing the boundedness of $(\partial \psi_3/\partial y)$ --
  and  (\ref{e_c14y}) -- its counterpart for $(\partial \psi_1/\partial y)$
   and $(\partial \psi_4/\partial y)$) --  it is easily derived that these equations are proportional   
  on the \cp and \cm curves if and only if for $t=t^\pm$
  $$  (t \gamo H + t \gamo E_c - \gam u \ww ) (u^2+ v^2)= ( (\gam+1) u^2 v - t \gamo u v^2 + (v + u t) ( \gamo E_c +\gamo H -\gam u^2))(u+vt)   $$ 
 Wherever $uv \neq 0 $, this equation is equivalent to
  $$  \gam (u^2 + v^2) =  ((\gam+1) u - t \gamo v)(u+vt) + ( \gamo E_c +\gamo H -\gam u^2)(1 + t^2)   $$
 that is found to be equivalent to  equation \colb{(\ref{e_deg2})}, the degree two equation
 $t$ which roots are
 the slope coefficients  of the \cp and \cm. Equation (\ref{e_c3yb}) hence also reduces to (\ref{e_c14y}) along the \cp and \cm curves.\\
  Finally, we consider the last differential form (\ref{e_c2yb}), expressing the boundedness of
  $(\partial \psi_2/\partial y)$. Whether it is proportional to (\ref{e_c14y}) is not straightforward in particular
  due to the complex expression of $C^2_{2y}$ and the ratio  $C^3_{2y}/C^2_{2y}$ that is not obviously equal to $t$.
  Nevertheless, we have proven that, along the characteristic curves the differential form expressing
  the boundedness of $(\partial \psi_2/\partial x)$ and  $(\partial \psi_2/\partial y)$ (equations (\ref{e_c2xb})
  and (\ref{e_c2yb})) are proportional by a minus $t$ factor. So we may use this property to derive a simpler
  expression of $C^2_{2y}$ or prove the proportionality of  (\ref{e_c2xb}) with (\ref{e_c14y}) along the \cp
  and \cm curves. Whatever the approach the condition for proportionality reads
$$   t~(\gamo H + \gamo E_c + \gam v \ww)(u^2 + v^2)  =  (  \gamo (u t + v) ( E_c + H) + u v^2 t  - \gamo u^2 v  -\gam v^3) (u+vt), $$ 
that is found to be equivalent to \colb{(\ref{e_deg2}) wherever $uv \neq0$}.\\
In a final formal calculation verification, it was checked that, on the \cp and \cm curves,  (\ref{e_c1yb}),
   is proportional to (\ref{e_c14y}) (and we already know that
 equations (\ref{e_c2yb}), (\ref{e_c3yb}), (\ref{e_c4yb}) are proportional to (\ref{e_c1yb})
 along these curves).
 
 \subsection{\colb{Main results and extension to 3D}}
\colb{  We have found two differential equations, (\ref{e_atraja})  and (\ref{e_atrajb}), valid
along the streamtraces \st for the adjoint system:
$$  E_c ~\frac{d\psi_1}{ds} +   H ( u~ \frac{d\psi_2}{ds} + v~ \frac{d\psi_3}{ds})  +   ~H^2  \frac{d \psi_4}{ds}  = 0   \qquad \quad
 \frac{d\psi_1}{ds} + u ~ \frac{d\psi_2}{ds} +  v ~ \frac{d\psi_3}{ds} +  E_c \frac{d\psi_4}{ds} = 0.  $$
 They are the counterpart of the constant total enthalpy and constant entropy properties for the flow. They may be combined (for example to derive again $\overline{U}.\nabla \psi_1- H \overline{U}.\nabla \psi_4 = 0  $) or their coefficients may be expressed differently using the well-known equations satisfied by steady inviscid flows along a streamtrace.}
 \colb{Let us finally note that their straightforward 3D extensions (with natural notations),
 \beas
  E_c ~\frac{d\psi_1}{ds} +   H ( u~ \frac{d\psi_2}{ds} + v~ \frac{d\psi_3}{ds}  +  w~\frac{d\psi_4}{ds})+ H^2 ~ \frac{d \psi_5}{ds}  = 0
  \\
 \frac{d\psi_1}{ds} + u ~ \frac{d\psi_2}{ds} +  v ~ \frac{d\psi_3}{ds} +   w ~ \frac{d\psi_4}{ds} + E_c~ \frac{d\psi_5}{ds} = 0 ,
 \eeas
 are valid. This is the subject of Appendix C. The demonstration used in 3D
  also provides another method for deriving the 2D equations.}
  \\
  \colb{
 We have found one differential equation for the \cp
     and one differential equation for the \cm, equation (\ref{edo_c14y}) with relevant value of $t^\pm$ for each curve \colb{(\ref{eteps})}
     }
\colb{     
\beas
    &&(u+v t^\pm)~\frac{d\psi_1}{ds} +  (u^2 + v^2) (\frac{d\psi_2}{ds} + t^\pm \frac{d\psi_3}{ds} ) + H  (u+v t^\pm)~\frac{d\psi_4}{ds} = 0. \\
    &&t^+ = \frac{uv + c\sqrt{u^2 + v^2 -c^2} }{u^2-c^2}   \quad \textrm{for a \cp} \qquad \qquad 
    t^- = \frac{uv - c\sqrt{u^2 + v^2 -c^2} }{u^2-c^2}   \quad \textrm{for a \cm} 
\eeas 
}
\colb{
     They are the counterpart of the differential forms satisfied by primitive
     flow variables along the \cp and \cm. In the simpler cases,
     where the classical angular relations (\ref{e:kpkm}) are
     valid, these relations may be used to express differently the coefficients. We do not expect these relations involving a 2D slope, $t$, to admit an extension in 3D.
}
%
%
%
\section{\label{sec:l1_3} \colb{Assessment of the adjoint ODEs}}
\subsection{\colb{Consistency with known flow perturbation mechanisms}}
%
\colb{ The adjoint vector is known to express the influence of a flow perturbation on the associated QoI.
 Although discrete and continuous adjoint methods are nowdays common tools for shape optimization, flow control and receptivity-sensitivity-stability
	analysis, adjoint vectors are not easily interpreted. The reason is that a single adjoint component, at a given location,
 is the rate of the change in the QoI to the amplitude of a local perturbation in the corresponding flow equation only (eg for the first component,
  mass injection without perturbation of the momentum and energy equations). These individual equation perturbations, of course, do not correspond to
 any realistic possibility. A second complementary point of view, already mentioned in the introduction,
 consists in calculating the dot product of the adjoint components with the vector of a realistic perturbation
 and discuss the map of the actuation influence \cite{GilPie_97,JPRenLab_22}.} 
\\
\colb{ Plots of individual adjoint components for Euler flows appear in a 2002 publication by Venditti and Darmofal \cite{VenDar_02} (fig. 19 and 26).
The discussion of a x-momentum $CL$-adjoint
 plot in \cite{VenDar_02}
  mentions {\sl a singularity in the adjoint along the stagnation streamline and a weak discontinuity upstream of the primal shock on the upper surface.}
(With a finer mesh, the latter would have been identified as a \cm impinging the upperside shock-foot).
	Although not discussed by the authors of  \cite{VenDar_02}, a corresponding plot for a supersonic flow about two airfoils,
	strongly suggest backward information
 propagation along the \st, \cm and \cp curves from the support of the QoI. 
  Concerning transonic airfoil flows, reference 
  \cite{TodVonBou_16}
  describes the mechanism by which locally perturbing 
  one component of the flow along the \cm (resp. \cp) impinging the upperside (resp. lowerside) shock-foot results in a
 strong change in the lift and drag: the flow perturbation propagates along the \cm (resp. \cp) and results in a displacement of the shock. 
}
\\
  \colb{ The evaluation of the influence of physical source terms on the lift
  or drag of a profile goes back 
  to Giles and Pierce \cite{GilPie_97} who introduced the four physical local source terms recalled in section I.
  Consistently with equations (\ref{e_atraja}), (\ref{e_atrajb}) and 
  (\ref{edo_c14y})
   the results obtained with this approach also support adjoint information propagation
  along the \st, \cm and \cp, from the support of the QoI (backwards
  w.r.t. the direction of flow perturbations). Regarding the specific goal of acting at a shockfoot and the four aforementioned source terms \cite{GilPie_97},
	the authors of \cite{JPRenLab_22} demonstrate that the  first two source terms (mass source at stagnation conditions, and normal force)
	may displace the shock and strongly alter near-field forces if located along the \cp / \cm of the shock-foot whereas the fourth source (change
	 in stagnation pressure at fixed static pressure and total enthalpy) may displace the shock-foot if located along the stagnation streamline or along the
	  wall upwind the shock.
}
\\	
\colb{The demonstrated ODEs along the \st, \cm and \cp curves are hence consistent with known lines of specific influence on drag or lift
 of classical steady Eulerian flows. These lines also appear in the search of optimal forcings in control studies \cite{Sar_14,SarMetSip_15} but we do not extend on this aspect
  due to the different base equations.}	 
  
\subsection{\colt{Consistency with the analytical adjoint field of 2D 
 supersonic constant flow areas and the equations for the adjoint gradient at shocks  }}
\begin {figure}[htbp]
  \begin{center}
	  \includegraphics[width=0.8\linewidth]{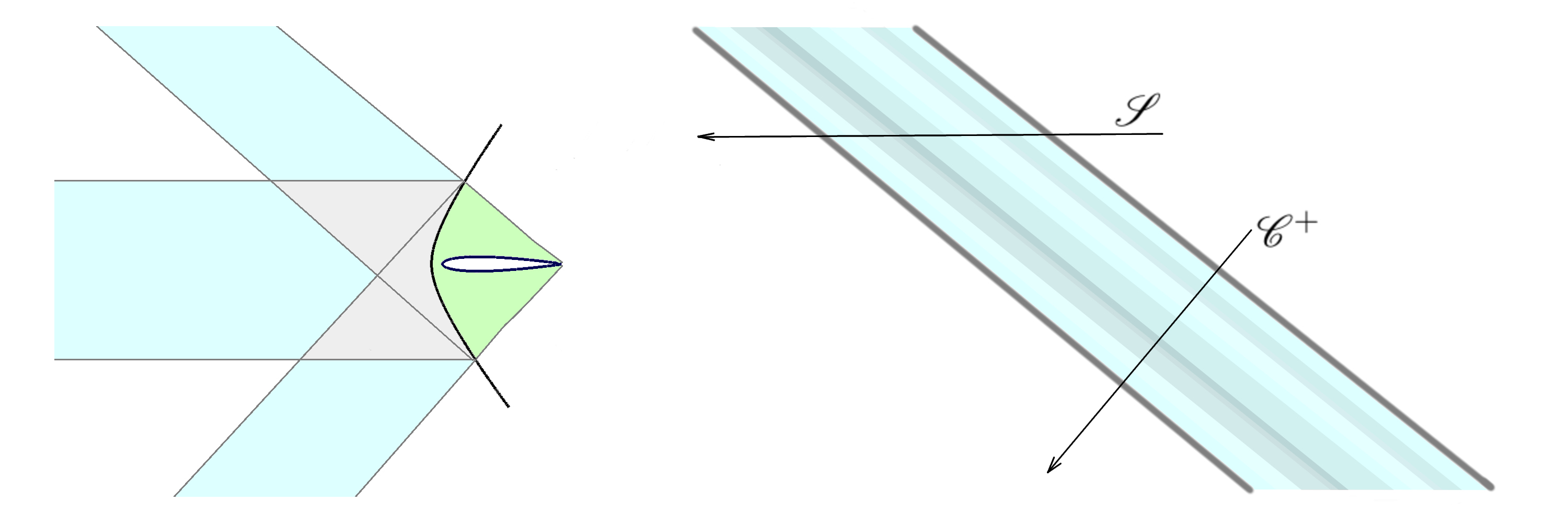}
	  \caption{\colb{Left: Sketch of lift/drag adjoint field for a supersonic
	   flow about an airfoil. Individual stripes (blue) and superimposition of the three stripes (grey) in the constant flow zone ; nonzero adjoint zone downwind the shockwave (green). Right:  \cm stripe crossed by a \st and a \cp curve}}
  \label{f:sksup}
  \end{center}
\end {figure}

%
\colb{In reference \cite{TodVonBou_16} Todarello {\sl et al.} derived the
  mathematical expression of the 2D Eulerian adjoint vector
 in a supersonic zone with constant flow (typically upwind the detached shockwave created by an airfoil). The angle of attack being $\alpha$ and
 Mach number $M$, this formula reads}
 \colb{
 \bea
   \psi(x,y) &=&  \varphi_{\alpha} (x ~\sin(\alpha) - y ~\cos(\alpha) ) \lambda^{\alpha}_0 +\varphi_{\alpha+\beta} (x ~\sin(\alpha+\beta) - y~ \cos(\alpha+\beta) ) \lambda^{\alpha+\beta}_0 \nonumber\\
   &+& \varphi_{\alpha-\beta} (x ~\sin(\alpha-\beta) - y ~\cos(\alpha-\beta) ) \lambda^{\alpha-\beta}_0,  \label{e:stripes} 
 \eea
 }
 \colb{where $\beta=\textrm{sin}^{-1}(1/M)$, $\varphi_{\alpha}$,
 $\varphi_{\alpha-\beta}$, $\varphi_{\alpha-\beta}$ are three scalar functions, the $\lambda_0^\mu$ are 
 left eigenvectors\cite{JPRenLab_22} of $A \sin(\mu) -B \cos(\mu)$}
{\footnotesize\colb{
 \beq
\lambda_0^{\alpha-\beta} = \begin{pmatrix}
   \frac{c}{\rho} \left( 1 + \frac{\gamma_1}{2} M^2 \right)\\
   \frac{1}{\rho} \left( \sin(\alpha-\beta) - \gamma_1\ds\frac{u}{c} \right)\\
   \frac{1}{\rho} \left( -\cos(\alpha-\beta) - \gamma_1\ds\frac{v}{c} \right) \\
   \frac{\gamma_1}{\rho c}
\end{pmatrix}, \quad
\lambda_0^{\alpha+\beta} = \begin{pmatrix}
   \frac{c}{\rho} \left( 1 + \frac{\gamma_1}{2}M^2 \right)\\
   -\frac{1}{\rho} \left( \sin(\alpha+\beta) + \gamma_1\ds\frac{u}{c} \right)\\
   -\frac{1}{\rho} \left( -\cos(\alpha+\beta) + \gamma_1\ds\frac{v}{c} \right)\\
   \frac{\gamma_1}{\rho c}
\end{pmatrix} \quad
\lambda_0^{\alpha} =  \begin{pmatrix}
   -1 - \frac{\gamma_1}{2}M^2\\
   \frac{\gamma_1 u}{c^2} + \frac{2~ \cos(\alpha)}{\cos(\alpha) u+\sin(\alpha)v}\\
   \frac{\gamma_1 v}{c^2} + \frac{2 ~\sin(\alpha)}{\cos(\alpha) u+\sin(\alpha)v}\\
   -\frac{\gamma_1}{c^2}.
\end{pmatrix}
    \label{e_lefteig}
    \eeq
    }}
    \colb{
which formulas \cite{Hir_07} have been simplified here 
using the null eigenvalues relations valid in this specific context:
$ u~\sin(\alpha)-v~\cos(\alpha)=0$, $\quad u~\sin(\alpha-\beta)-v~\cos(\alpha-\beta)+c =0 $, $\quad u~\sin(\alpha+\beta)-v~\cos(\alpha+\beta)-c =0$. }
\colb{Equation (\ref{e:stripes}) is the mathematical formula for the three stripes field depicted in figure \ref{f:sksup}. 
Each $\varphi$ scalar function
 expresses the combined variation in the amplitude of the 
 adjoint components normally to the stripe direction. Each stripe is
  crossed by the characteristic lines oriented in the direction of the other two and we 
  question whether equation (\ref{e_atraja}), (\ref{e_atrajb}), (\ref{edo_c14y}) provide new
   information on the $\varphi$ functions.}
  \\
  \colt{Considering the stripe oriented in the $\alpha-\beta$ direction,
 we first note that equation   (\ref{edo_c14y}) with the $t^-$ value,
 is automatically 
  satisfied in its geometrical domain since 
  $\varphi_{\alpha-\beta} (x ~\sin(\alpha-\beta) - y ~\cos(\alpha-\beta) ) \lambda^{\alpha-\beta}_0$ induces no variation
   of the $\psi$ components in the \cm direction.
  Regarding the conditions for satisfying (\ref{e_atraja}) and (\ref{e_atrajb})
   along an \st curve, and satisfying (\ref{edo_c14y}) 
   along a \cp curve where they cross the $\alpha-\beta$ stripe
  (as in fig. \ref{f:sksup} right), we introduce the three functions 
}
\colt{
\beas
\Gamma^1_{S}(s)&=& E_c \psi_1 + H u \psi_2 + Hv \psi_3 + H^2 \psi_4\\
\Gamma^2_{S}(s)&=& \psi_1 + u \psi_2 + v \psi_3 + E_c \psi_4\\
\Gamma_{C+}(s)&=& (u+v t^+)\psi_1 + (u^2 + v^2) (\psi_2 + t^+  \psi_3) + H (u+ v t^+)  \psi_4,
\eeas
}
\colt{with $s$ the curvilinear abscissa along the curve mentioned in the index. It is proven in
appendix D that
$$ \frac{d\Gamma^1_{S}}{ds} = 0 \qquad \frac{d\Gamma^2_{S}}{ds} = 0  \qquad \frac{d\Gamma_{C+}}{ds} = 0 $$ without condition on the $\varphi$ functions. As the flow is constant, this is equivalent to the satisfaction of the ODEs along the \st~ and the \cp curves.
The other possible crossing of the stripes and the curves 
 also do not provide conditions on the $\varphi$ functions.
 We hence have not gained information on the
analytical adjoint field of 2D  supersonic constant flow zones
 but this automatic consistency establishes a link,
 via a series the orthogonality properties,
 between the
 coefficients of the differential equation  (\ref{e_atraja}), (\ref{e_atrajb}), (\ref{edo_c14y})
  and the relevant left eigenvectors of the Euler flux Jacobian
  (\ref{e_lefteig}). }\\

\colt{
  Regarding shock-waves, it is know that, for 
 classical QoIs like pressure integrals at
 the wall, the adjoint vector 
  is continuous at shocks although the flow is not, whereas
 the adjoint gradient may be discontinuous. In most common cases,
  although a shockwave only makes the normal
component of the velocity subsonic, the \cp and \cm curves
end at the shock and the most interesting discussion regards the 
 consequences of (\ref{e_atraja}), (\ref{e_atrajb}). As they are valid
  both sides of the shock, the jump operator may be applied to them
  across the discontinuity $\Sigma$:
        $$ \llbracket E_c ~\frac{d\psi_1}{ds} +   H ( u~ \frac{d\psi_2}{ds} + v~ \frac{d\psi_3}{ds})  +   ~H^2  \frac{d \psi_4}{ds} \rrbracket = 0   \qquad
 \llbracket \frac{d\psi_1}{ds} + u ~ \frac{d\psi_2}{ds} +  v ~ \frac{d\psi_3}{ds} +  E_c \frac{d\psi_4}{ds} \rrbracket = 0,   $$
 where all terms but $H$ may be discontinuous. Unfortunately, it 
 does not seem possible
 to reduce these equations to the simpler ones 
 for the adjoint gradient discontinuity \cite{Loz_18,JPRenLab_22} where the derivatives 
 in the two directions of the local frame of reference attached to
  the shock appear independently.
  Conversely, we note that in the case of a normal shock, the
  equations derived in \cite{Loz_18} prove the previous two
   jump equations.
}
%
\subsection{Numerical assessment method}
The numerical assessment method consists in computing flow and adjoint fields
over a very fine mesh, and calculating the following integrals:
\bea
 K_OS^1 &=& \int_{\cal S} \left(   E_c ~\frac{d\psi_1}{ds} +   H ( u~ \frac{d\psi_2}{ds} + v~ \frac{d\psi_3}{ds})  +   ~H^2  \frac{d \psi_4}{ds}  \right) ds    \\
 K_OS^2 &=& \int_{\cal S} \left( \frac{d\psi_1}{ds} +    u~ \frac{d\psi_2}{ds} + v~ \frac{d\psi_3}{ds}  +   ~E_c  \frac{d \psi_4}{ds} \right)  ds  \\
 K_OC^+ &=& \int_{{\cal C}^+}  \left(  (u+v t^+)~\frac{d\psi_1}{ds} +  (u^2 + v^2) (\frac{d\psi_2}{ds} + t^+ \frac{d\psi_3}{ds} ) + H  (u+v t^+)~\frac{d\psi_4}{ds}  \right) ds\\
 K_OC^- &=& \int_{{\cal C}^-}  \left(  (u+v t^-)~\frac{d\psi_1}{ds} +  (u^2 + v^2) (\frac{d\psi_2}{ds} + t^- \frac{d\psi_3}{ds} ) + H  (u+v t^-)~\frac{d\psi_4}{ds} \right) ds.
\eea
Here the intermediate subscript $O$ stands for the output functional of interest ;
it is subsequently replaced by $L$ (for the lift, $CLp$) and $D$
 (for drag, $CDp$,
 consistently with the adjoint vector placed on the right-hand side.

\begin {figure}[htbp]
  \begin{center}
	  \includegraphics[width=0.32\linewidth]{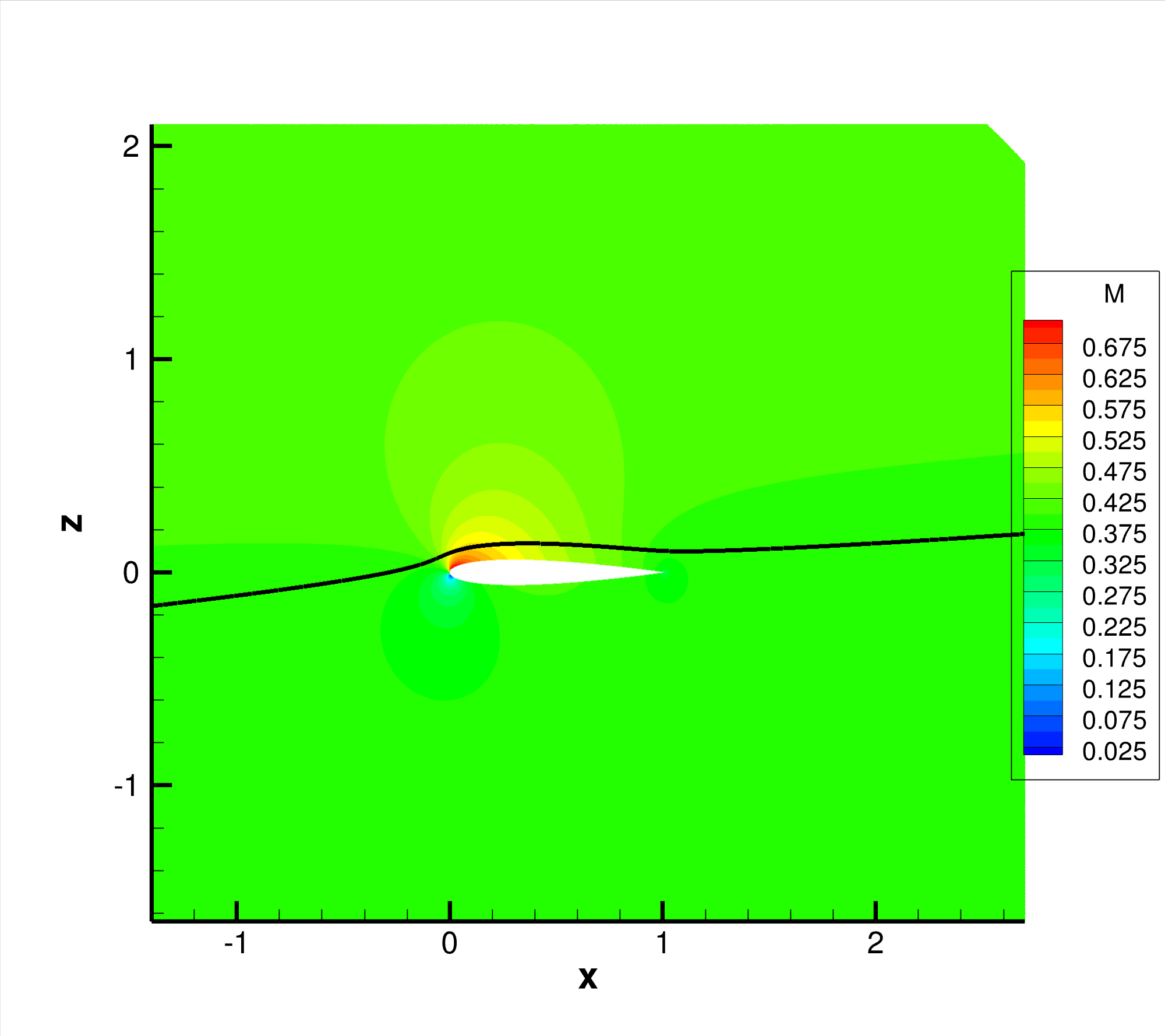}
	  \includegraphics[width=0.32\linewidth]{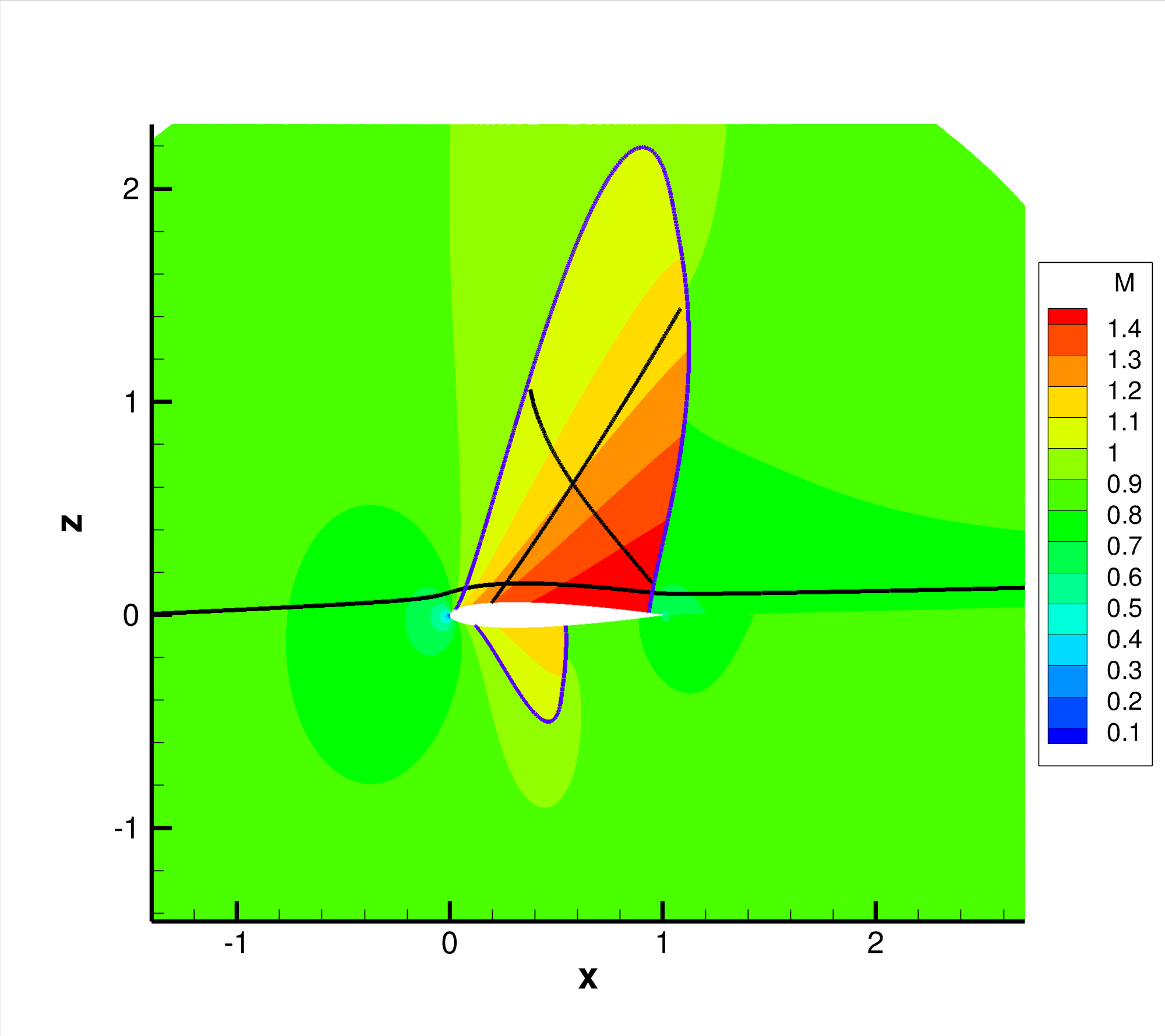}
	  \includegraphics[width=0.32\linewidth]{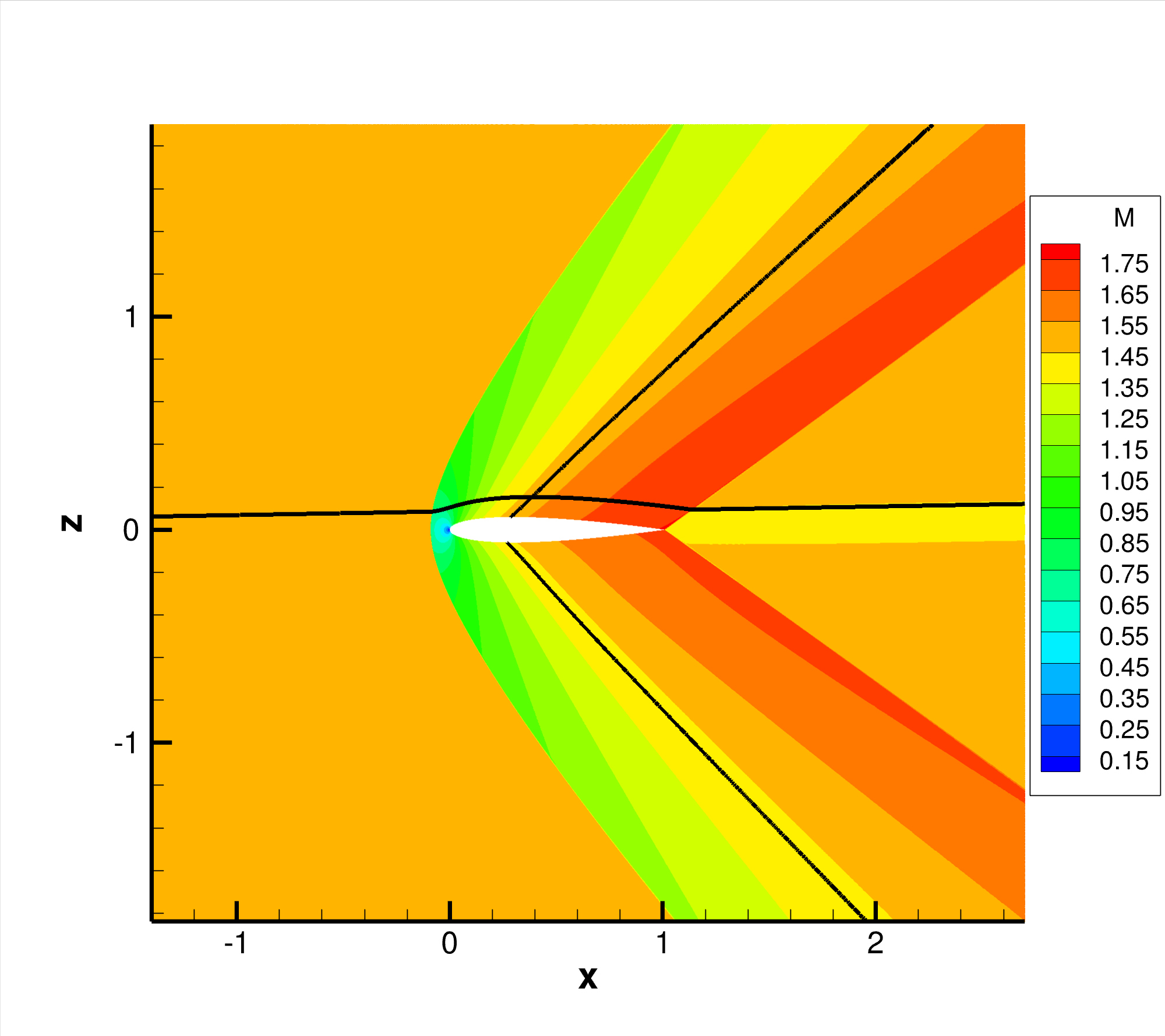}\\
	  \includegraphics[width=0.32\linewidth]{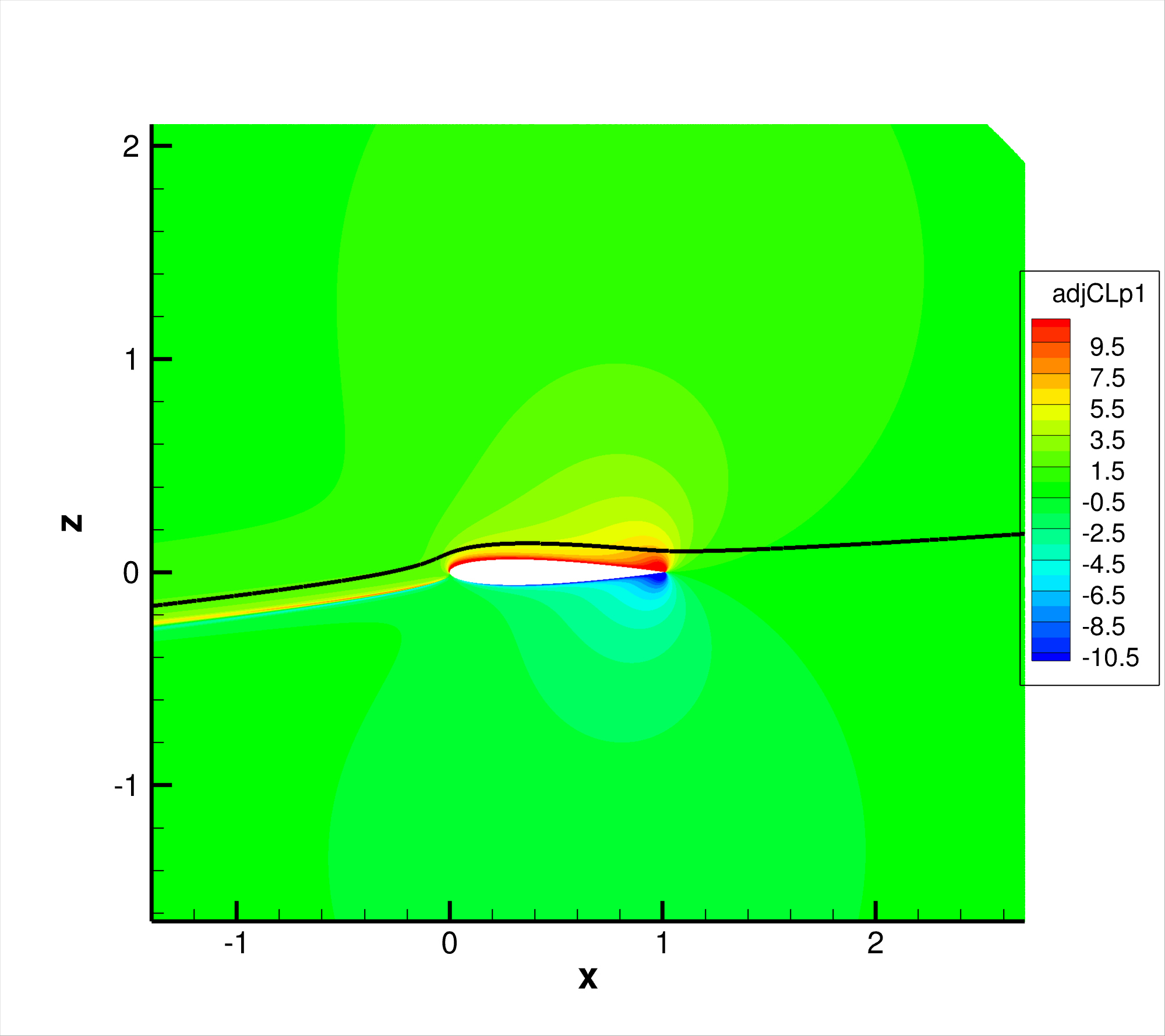}
	  \includegraphics[width=0.32\linewidth]{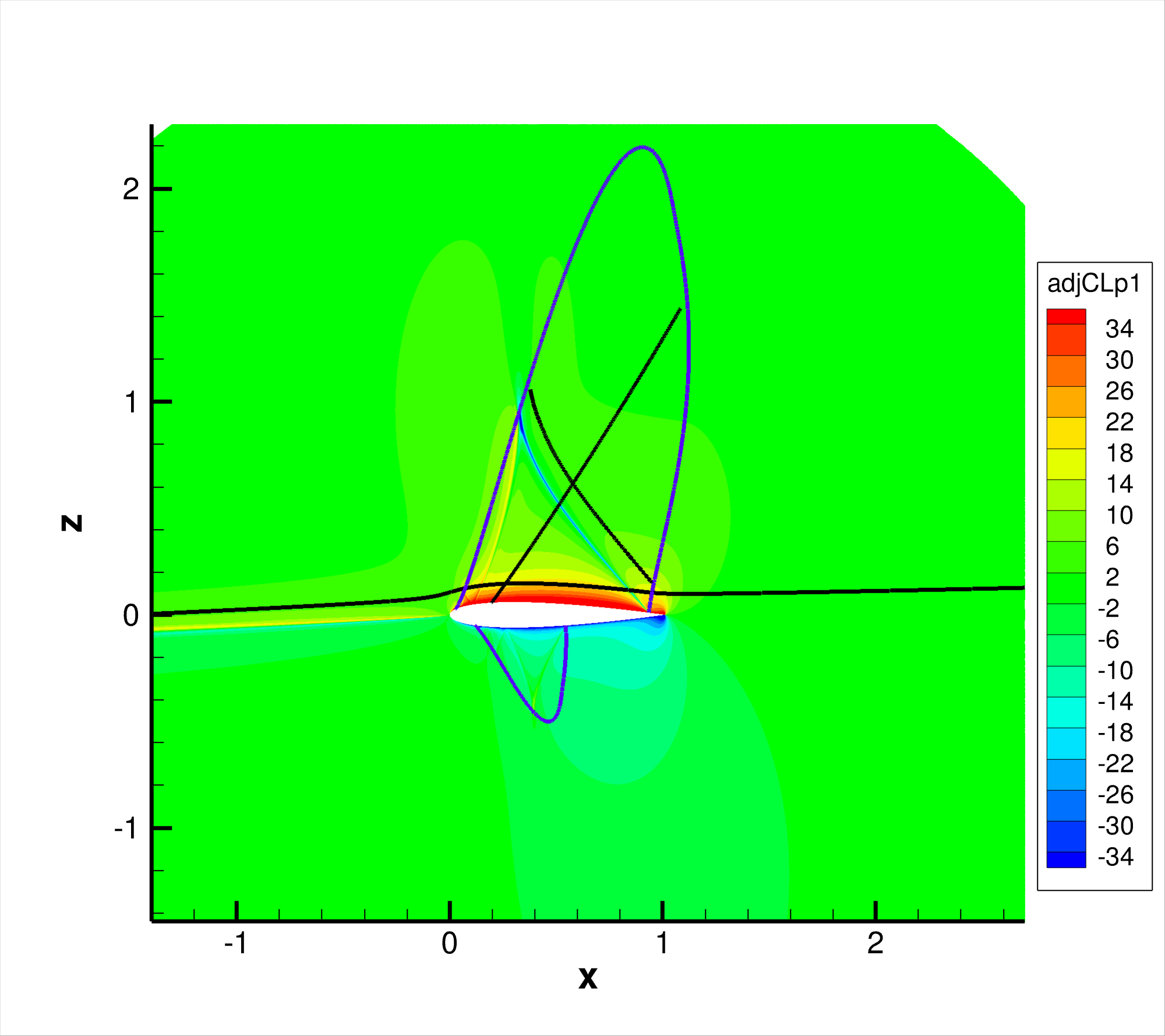}
	  \includegraphics[width=0.32\linewidth]{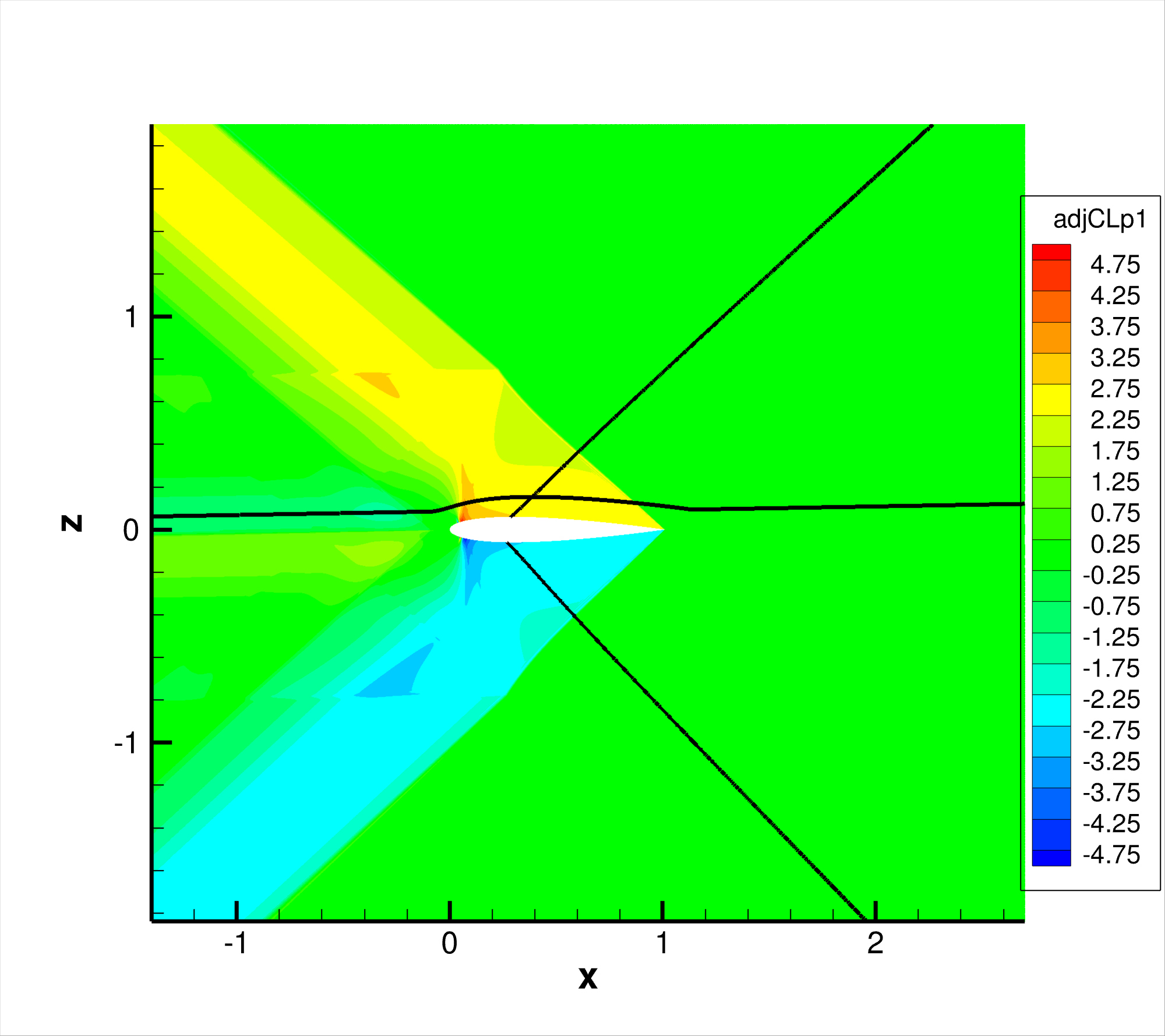}\\
	  \caption{Simulations over the 4097$\times$4097 mesh of \cite{VasJam_10}. iso-Mach number lines \colb{(upper three plots) and iso-$\psi^1_{CLp}$ lines 
	  (lower three plots)} and extracted
      curves. Left, subsonic case:  streamline. Middle, transonic case:  streamline, \cp, \cm (black) and sonic line (violet).
      Right supersonic case: streamline, \cp (upper side), \cm (lower side) }
  \label{f:str_mcurves}
  \end{center}
\end {figure}

The integration is performed  in the forward sense for the adjoint, that is,
backwards w.r.t. the direction of the flow information propagation. 
The integration domain for the above line integrals extends to the interior of the disk of radius $3c$ 
centred at $(0.5c,0)$, chosen for plotting readability, while the flow computational domain itself 
extends to 150$c$.
  It may be shorter, in particular in the transonic
 case where the \cp and \cm curves are limited to the supersonic bubble(s). 
 The four quantities $K_OS^1$,..., $K_OC^-$ are expected to be close to zero and, 
to avoid any error in scale, also calculated and plotted are the corresponding subparts, that is, for $ K_OS^1$ for example, 
$$ K_OS^1_1 = \int_{\cal S}  \left( E_c ~\frac{d\psi_1}{ds} \right) ds \qquad
K_OS^1_2 = \int_{\cal S} \left(  H u~ \frac{d\psi_2}{ds}  \right) ds $$
$$ K_OS^1_3 = \int_{\cal S} \left( H v~ \frac{d\psi_3}{ds}  \right) ds \qquad
K_OS^1_4 = \int_{\cal S} \left( H^2  \frac{d \psi_4}{ds} \right) ds .  $$
The sum of the four terms is expected to be much smaller than each one of them individually.
All the integrals are calculated backwards, along a finely discretized characteristic curve,
  simply by the trapezoidal rule.\\   
 The discrete flows \colt{and adjoints} were available from former computations \cite{JPRenLab_22} in which
the Jameson-Schmidt-Turkel scheme \cite{JamSchTur_81} was applied, and 
using the discrete adjoint module of the $elsA$ code \cite{JPDwi_10,JPRenDum_15}. 
\colt{Of course when trying to assess properties of exact adjoint
fields from numerical discrete solutions, it is desirable to work
  either with continuous adjoints or with dual consistent discrete adjoints\cite{MohPir_01,AlaPir_12,Loz_16}. Precisely in \cite{JPRenLab_22},}  it was demonstrated for structured meshes
 how to slightly modify the scheme's Jacobian (in the derivative of the dissipation flux, for the
next to wall faces) to get a dual consistent linearization.
This slight modification of the exact scheme Jacobian is retained here to work with adjoint fields that are
 consistent with the continuous equations discussed in \S 2. \colt{Note
 also that these adjoint fields have also been satisfactorily verified
  by a posteriori discretization of the continuous adjoint equation\cite{JPRenLab_22}}.\\
 Only the solutions calculated over the finest mesh defined in reference
 \cite{VasJam_10} (structured 4097$\times$4097 mesh) are used here. 
 The iso-Mach number lines, \colt{iso-first component of $CLp$ adjoint} and the extracted curves
    may be seen for all cases in figure \ref{f:str_mcurves}.
%
\subsection{\colb{Numerical assessment of the ODEs for a} supersonic flow about the NACA0012 airfoil}
\begin {figure}[htbp]
  \begin{center}
	  \includegraphics[width=0.4\linewidth]{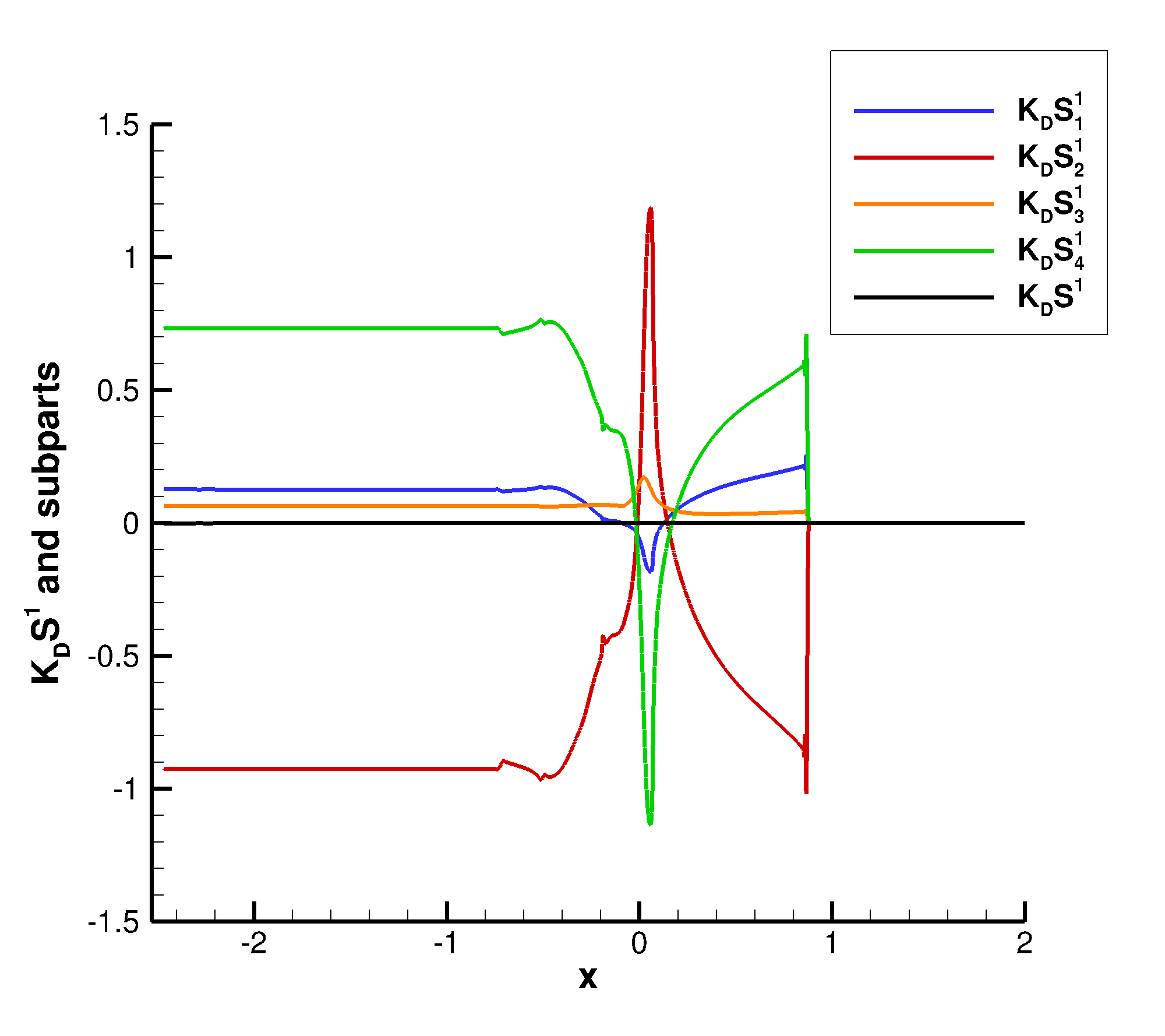}
	  \includegraphics[width=0.4\linewidth]{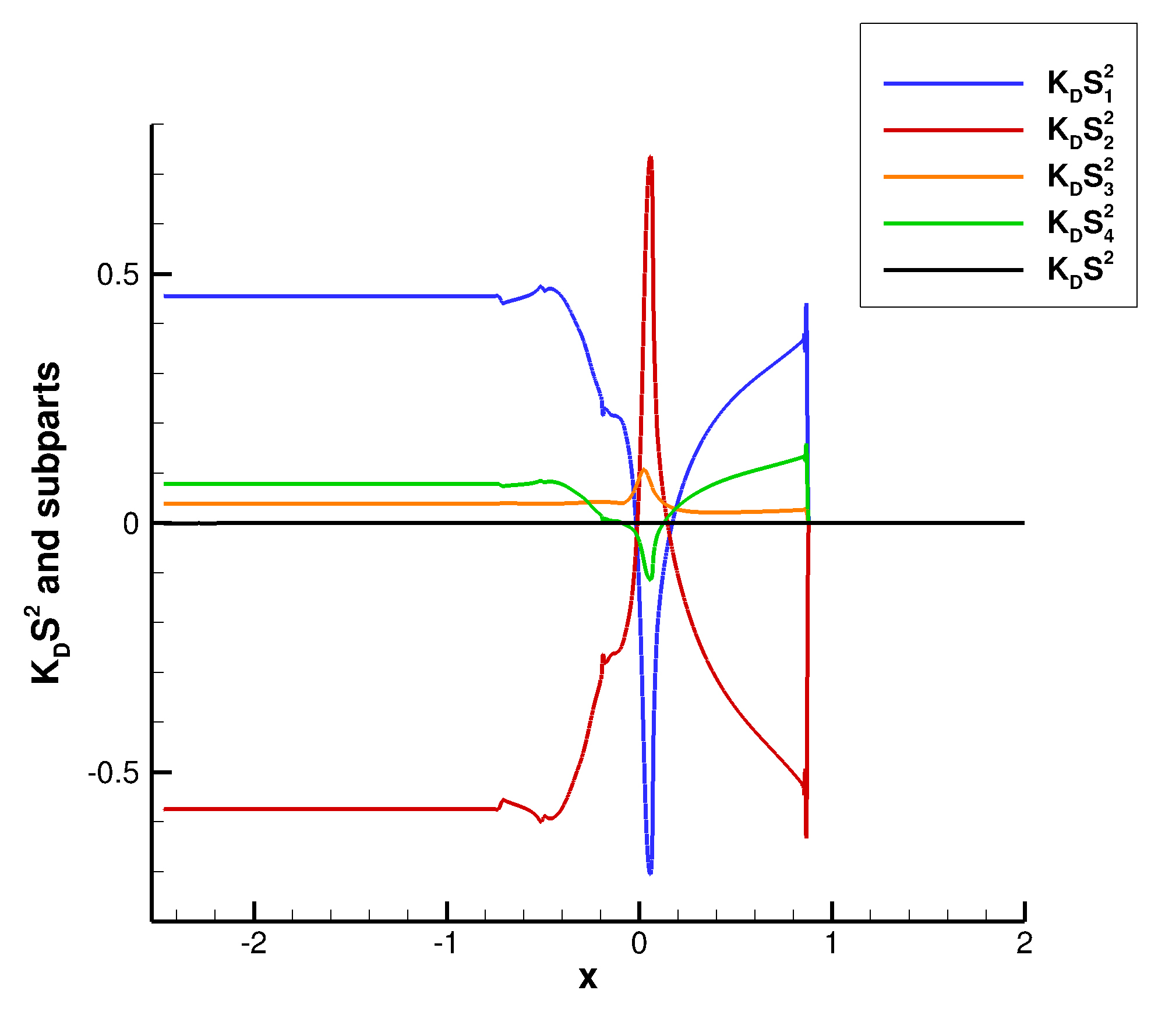}
	  \caption{$M_\infty=1.50$, $\alpha=1^o$, (4097$\times$4097 mesh \cite{VasJam_10}) Numerical assessment of equation (\ref{e_atraja}) (left)
      and (\ref{e_atrajb}) (right) for the lift.
      Method of verification: the black curve should ideally coincide with the $x$ axis }
  \label{f:str_sup}
  \end{center}
\end {figure}

\begin {figure}[htbp]
  \begin{center}
	  \includegraphics[width=0.4\linewidth]{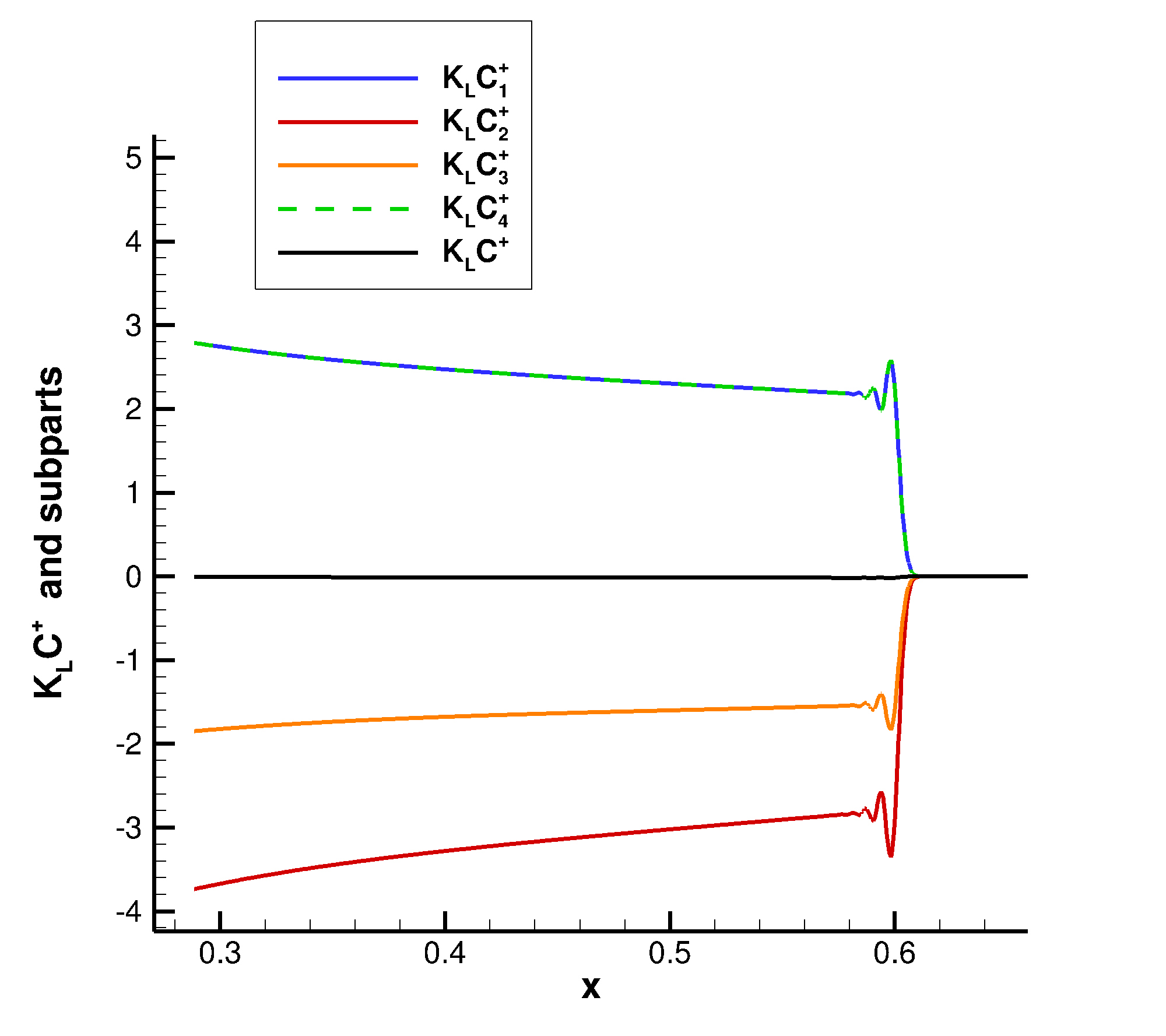}
	  \includegraphics[width=0.4\linewidth]{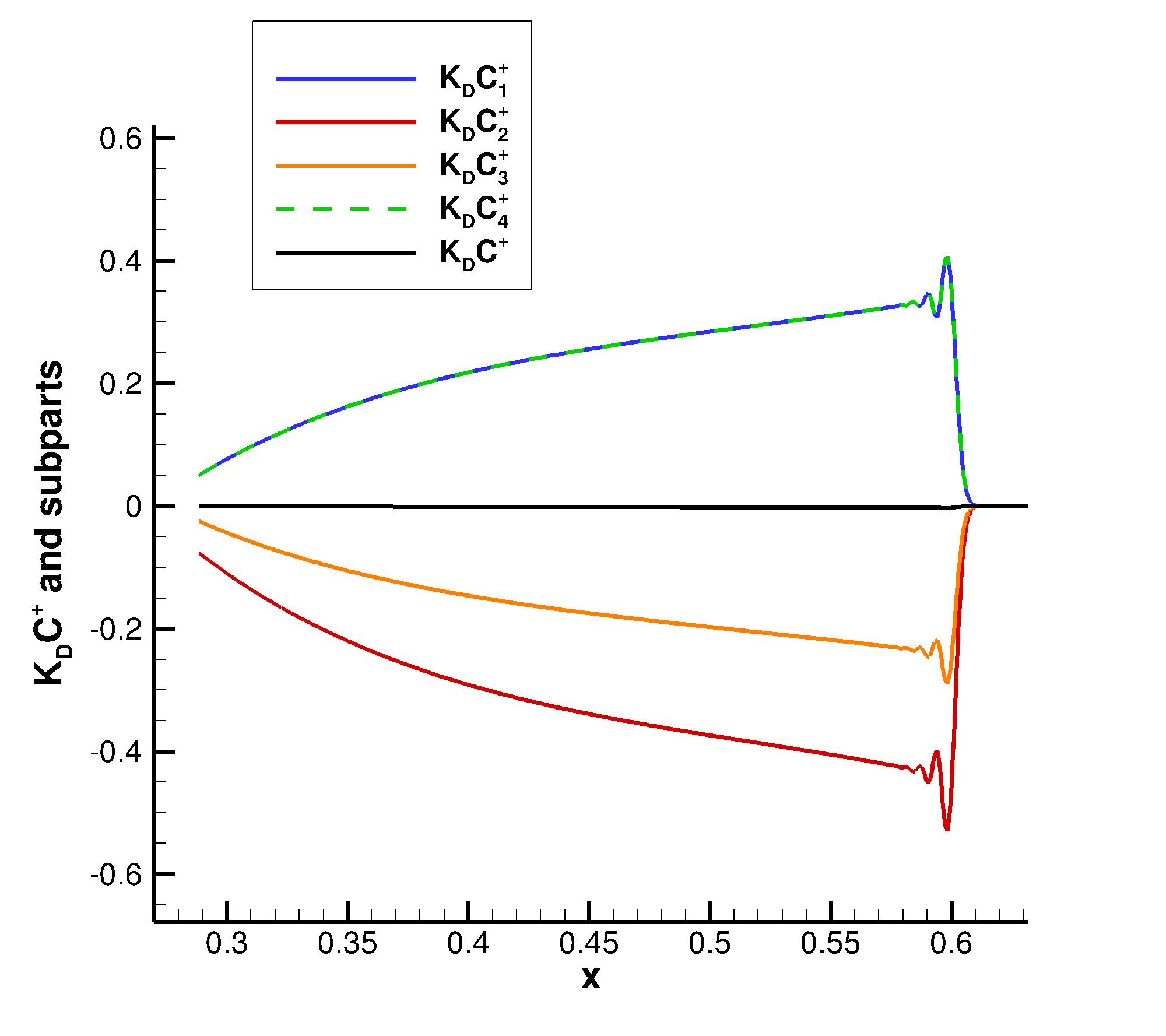}
	  \caption{$M_\infty=1.50$, $\alpha=1^o$, (4097$\times$4097 mesh \cite{VasJam_10})
      Numerical assessment of equation (\ref{e_c14y}) for the lift  (left) and for the drag (right).
      Method of verification: the  black curve should ideally coincide with the $x$ axis }
  \label{f:cpcm_sup}
  \end{center}
\end {figure}

 The retained flow conditions are $M_\infty=1.5~,~ \alpha=1^o$. 
We first assess the streamtraces equations (\ref{e_atraja}) (\ref{e_atrajb}). The  $K_DS^1$,
$K_DS^2$,  $K_LS^1$,  $K_LS^2$ integrals and their subparts are
calculated along the trajectory passing through $(c,0.1c)$. The integration indeed leads to very small
values of $K_DS^1$,
$K_DS^2$,  $K_LS^1$,  $K_LS^2$ along the curve w.r.t. their subparts.
It is well-known for this of kind of flow that the exact lift- and drag-
adjoint is equal to zero downstream the backwards flow-characteristics emanating from the trailing edge
(since no perturbation downstream those two lines can affect the pressure on the aerofoil and, consequently,
the lift or the drag -- see for example fig. 6 and A21 in \cite{JPRenLab_22}).
This property is well satisfied
by discrete adjoint fields and, as the integration
is performed backwards along the streamtrace,
null values of  $K_DS^1$, $K_DS^2$,  $K_LS^1$,  $K_LS^2$ and all their subparts are observed
above a specific $x$ corresponding to the intersection of the streamtrace with this
trailing-edge \cm. \colt{The integration of (\ref{e_atraja}) and (\ref{e_atrajb}) 
 reveals (i) the discontinuity of the 
integration variables (the adjoint components) at $x\simeq0.85$ that appears as a discontinuity in the subpart curves ; (ii)
a discontinuity of the integration coefficients, when crossing the detached shock wave at
 $x\simeq -0.08$, that results in strong gradients in the subparts curves.
  None of those discontinuities alters the almost null values of $K_DS^1$
  $K_DS^2$, $K_DL^1$ and $K_LS^2$.}
 \colt{The strong adjoint gradients in the subsonic bubble (approximately $x \in ~[-0.08,0.05]$)
  also clearly translate the $K_S$ curves. Finally, note that for $x$ 
   lower than
  $\simeq -0.05$ the flow is constant (so the streamtrace is a straight line) and 
   for $x$ lower than  $\simeq -1$ the backward streamtrace enters a zone of
   constant adjoint along its direction (see fig. \ref{f:sksup} left and equation (\ref{e:stripes})). This theoretical property is well 
   translated in constant sections of the curves.
 The results are equivalently accurate for lift and drag.
 They are presented for the drag in figure \ref{f:str_sup}.}\\
A \cp and a \cm curve are then extracted using equation (\ref{e:kpkm}). The selected
\cp is initiated at $x\simeq0.3$
upper side and the retained \cm starts at the same abscissa but on the lower side.
This choice was guided by the extraction method and the observation that $k^+$ (resp. $k^-$) 
is almost constant on the lower (resp. upper) side.
 The  $K_DC^+$, $K_LC^+$, $K_DC^-$,$K_LC^-$ terms and their subparts have been
computed. The results appear to be very satisfactory. Also observed is the equality
of $K_DC^+_1$ and  $K_DC^+_4$, $K_LC^+_1$ and  $K_LC^+_4$ and along the \cp and correspondingly 
along the \cm curves. This is due to the fact that $\psi_1=H\psi_4$ (for Euler flows and for
 pressure-based integrals along the wall that is well satisfied at
the discrete level -- see for example \cite{JPRenLab_22} fig. A21, A22, A23) and to the expression
of the $d\psi_1$ and $d\psi_4$ terms in (\ref{e_c14y}). Figure \ref{f:cpcm_sup} presents
 the verification plots for the two functions along the selected \cp.

\subsection{\colb{Numerical assessment of the ODEs for a} transonic flow about the NACA0012 airfoil}
\begin {figure}[htbp]
  \begin{center}
	  \includegraphics[width=0.4\linewidth]{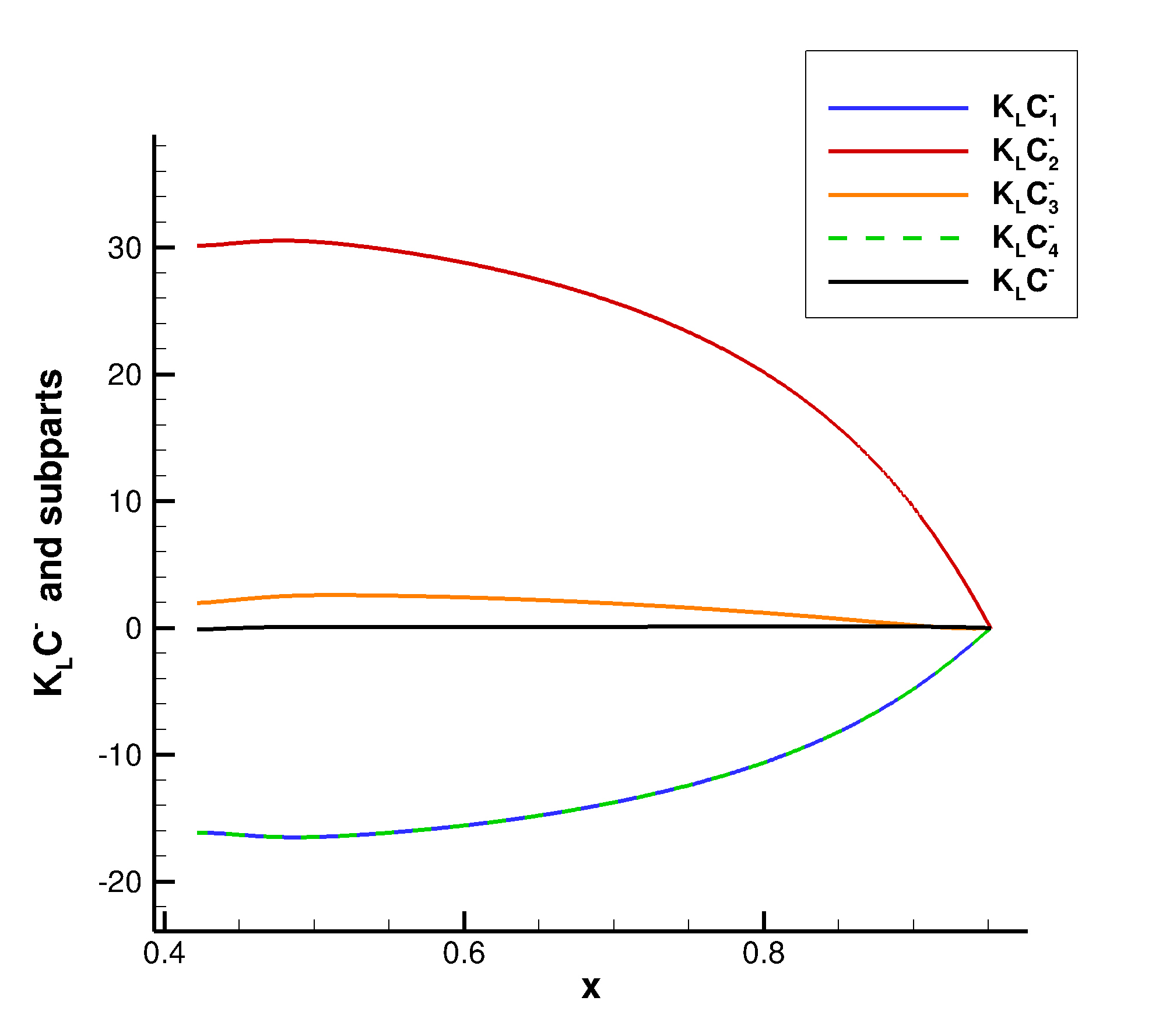}
	  \includegraphics[width=0.4\linewidth]{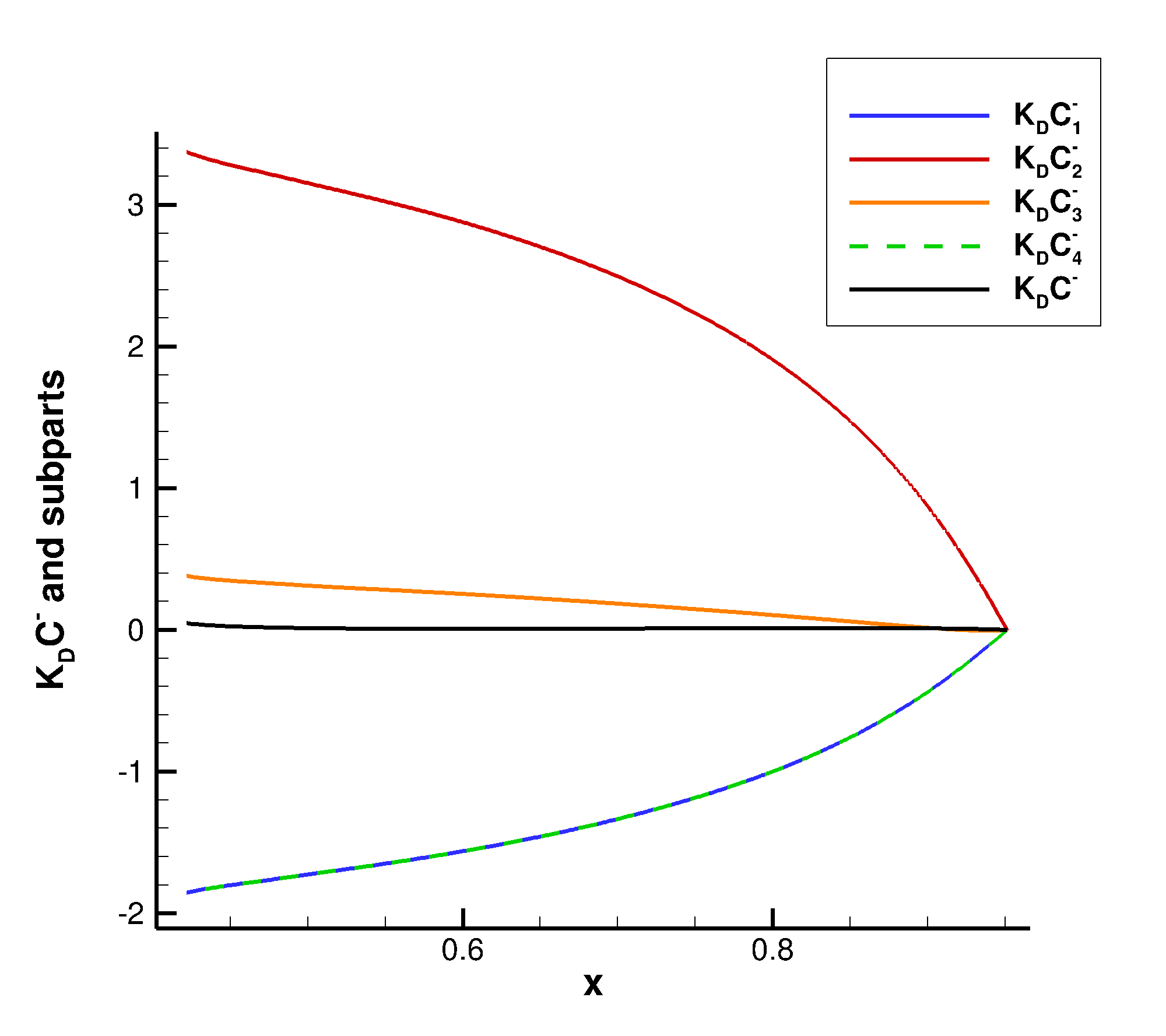}\\
    \includegraphics[width=0.4\linewidth]{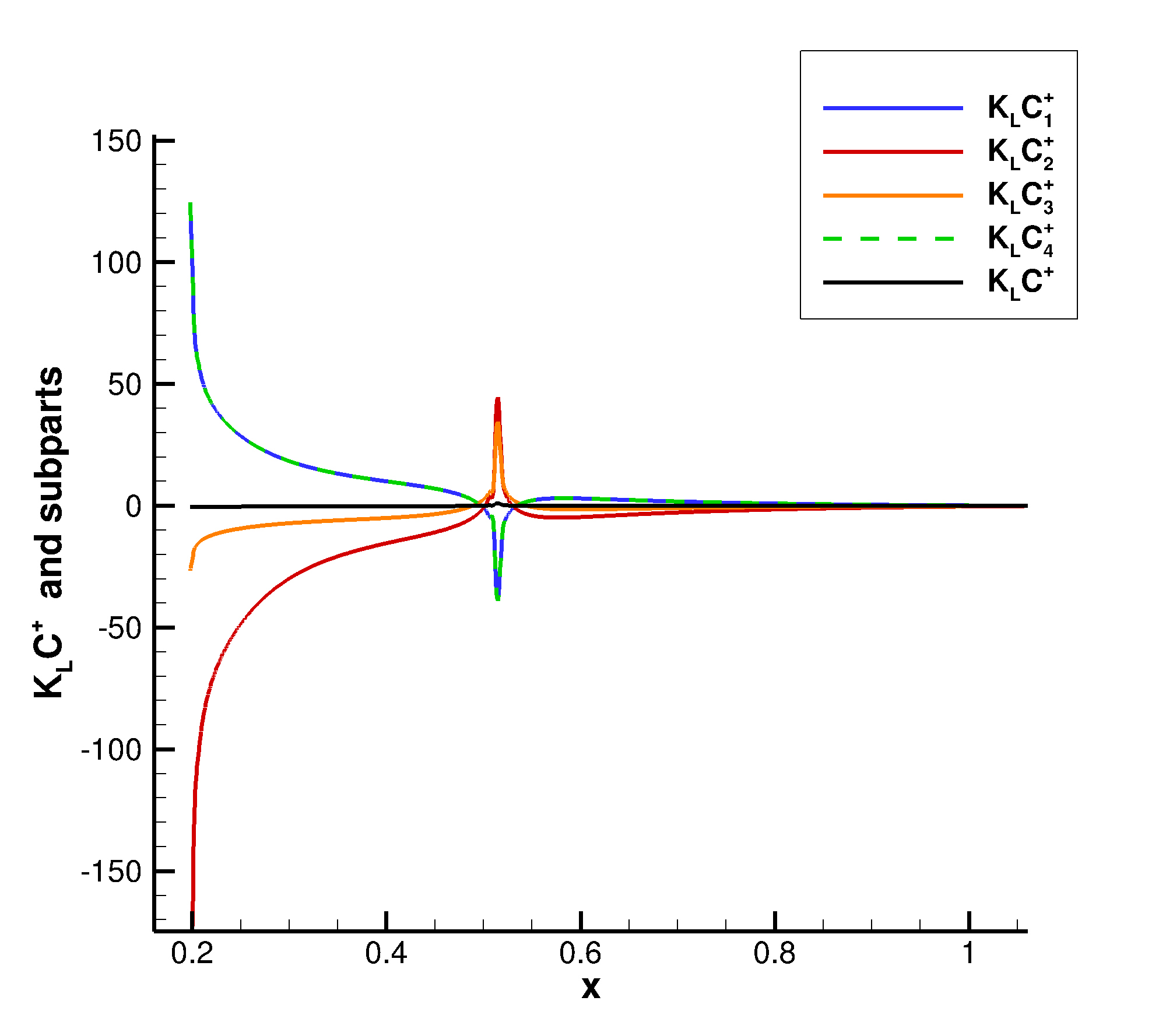}
	  \includegraphics[width=0.4\linewidth]{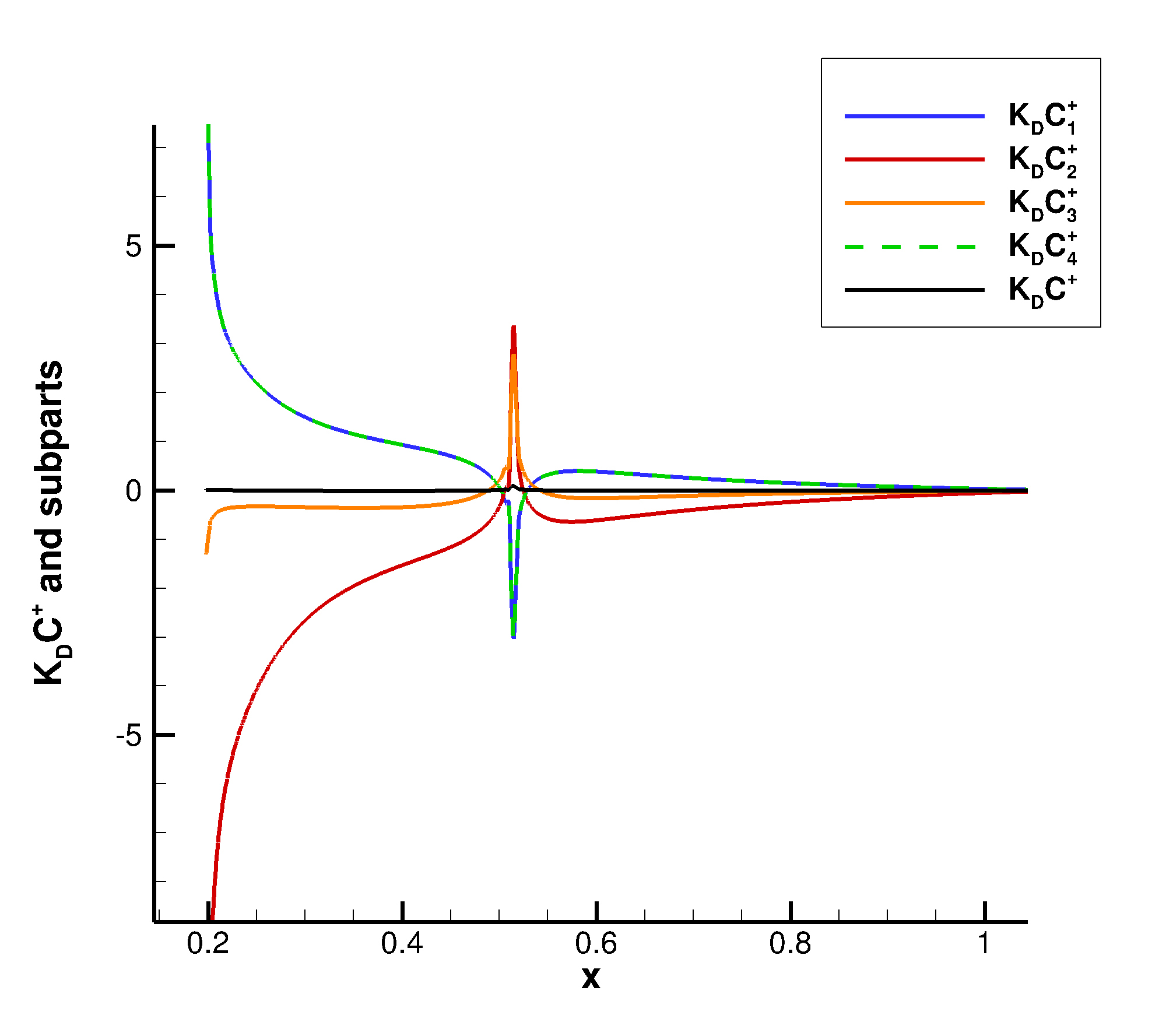}\\
	  \caption{\colb{$M_\infty=0.85$, $\alpha=2^o$}, (4097$\times$4097 mesh \cite{VasJam_10})
      Numerical assessment of equation (\ref{e_c14y}) for the lift  (left) and for the drag (right), \colb{for the selected \cm (up) and \cm (down)}.
      Method of verification: the  black curve should ideally coincide with the $x$ axis }
  \label{f:cpcm_trans}
  \end{center}
\end {figure}

The flow conditions are $M_\infty=0.85~,~ \alpha=2^o$.
Careful verification of the streamstraces equations (\ref{e_atraja}) (\ref{e_atrajb})
has been performed for the streamtrace passing through
 $(c,0.1c)$ and very satisfactory results have been found. \colb{As in the previous section,
  the intersection of the \st curve with the shockwave results in the 
   subsparts curves in very strong gradients but does not affect the almost null
    value of  $K_DS^1$, $K_DS^2$,  $K_LS^1$,  $K_LS^2$.}
As similar results have been presented in the previous and the next subsections, we focus here on
the \cp and \cm curves. A  \cp and a \cm curves have been extracted
taking care to select the longest
possible curves (and to avoid, for the \cm the curve passing by the shock-foot where numerical
divergence of the adjoint may be observed). The selected \cp (resp. \cm) is passing approximately
  through the point  $(0.197,0.057)$ (resp. $(0.954,0141)$). 
  The verification of the consistency of the numerical solutions w.r.t.
 (\ref{e_c14y}) is satisfactory although the largest observed errors
  appear in this case, for the \cm curve, \colt{for the lift}, in the vicinity of the inlet of the supersonic bubble 
-- see figure \ref{f:cpcm_trans}. This largest observed error
 is about 2\% of the  largest absolute value of the four subparts. 
 The integrals along the \cm are regular, whereas those along the 
 \cp exhibit a sharp peak close to $x \simeq 0.52$, at the intersection with
  the \cm passing by the shockfoot. We refer to \S3.1 for the reason of the
   corresponding strong adjoint gradient.

%
\subsection{\colb{Numerical assessment of the streamtrace ODEs for a} subsonic flow about the NACA0012 airfoil}
\begin {figure}[htbp]
  \begin{center}
	  \includegraphics[width=0.4\linewidth]{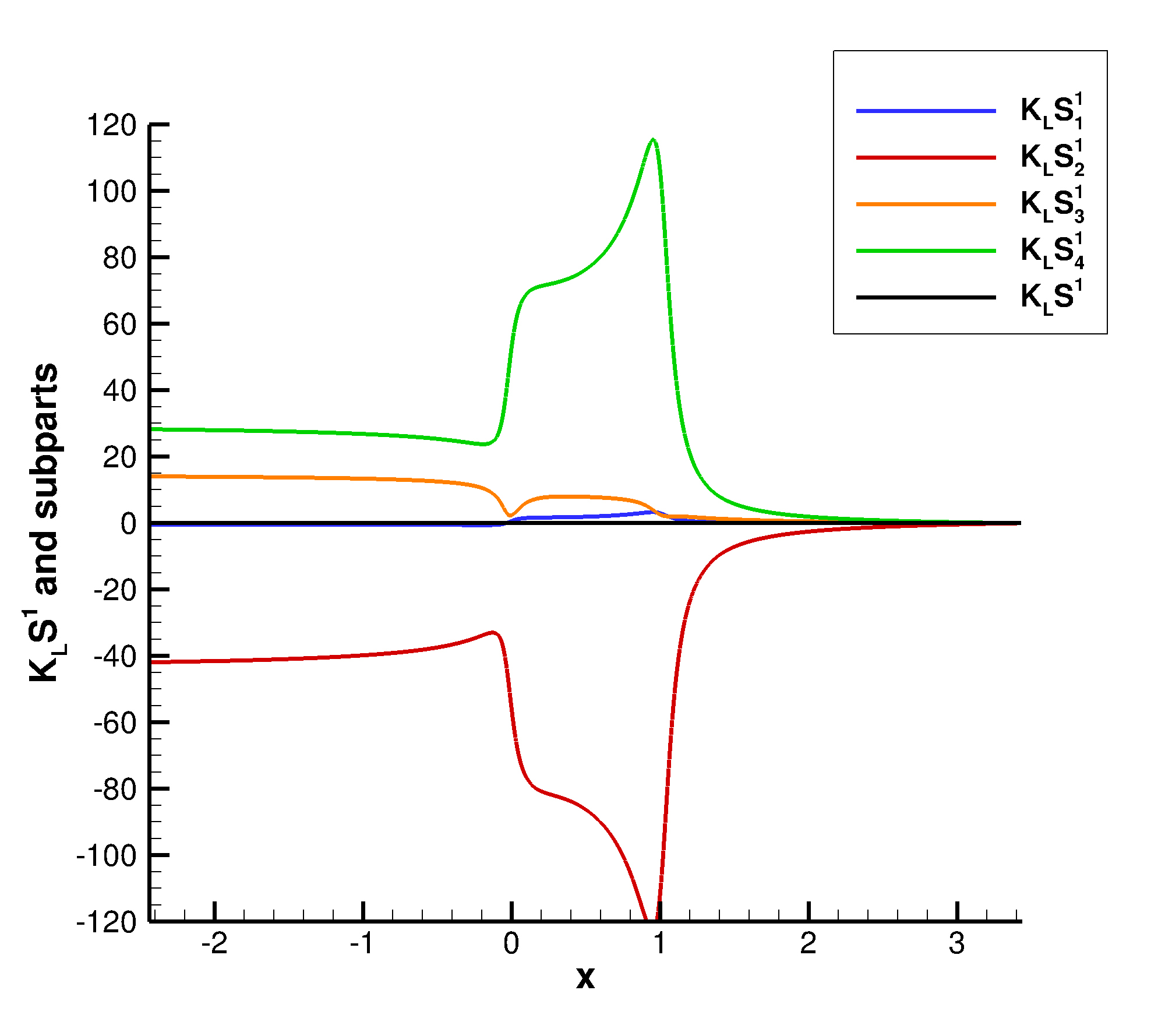}
	  \includegraphics[width=0.4\linewidth]{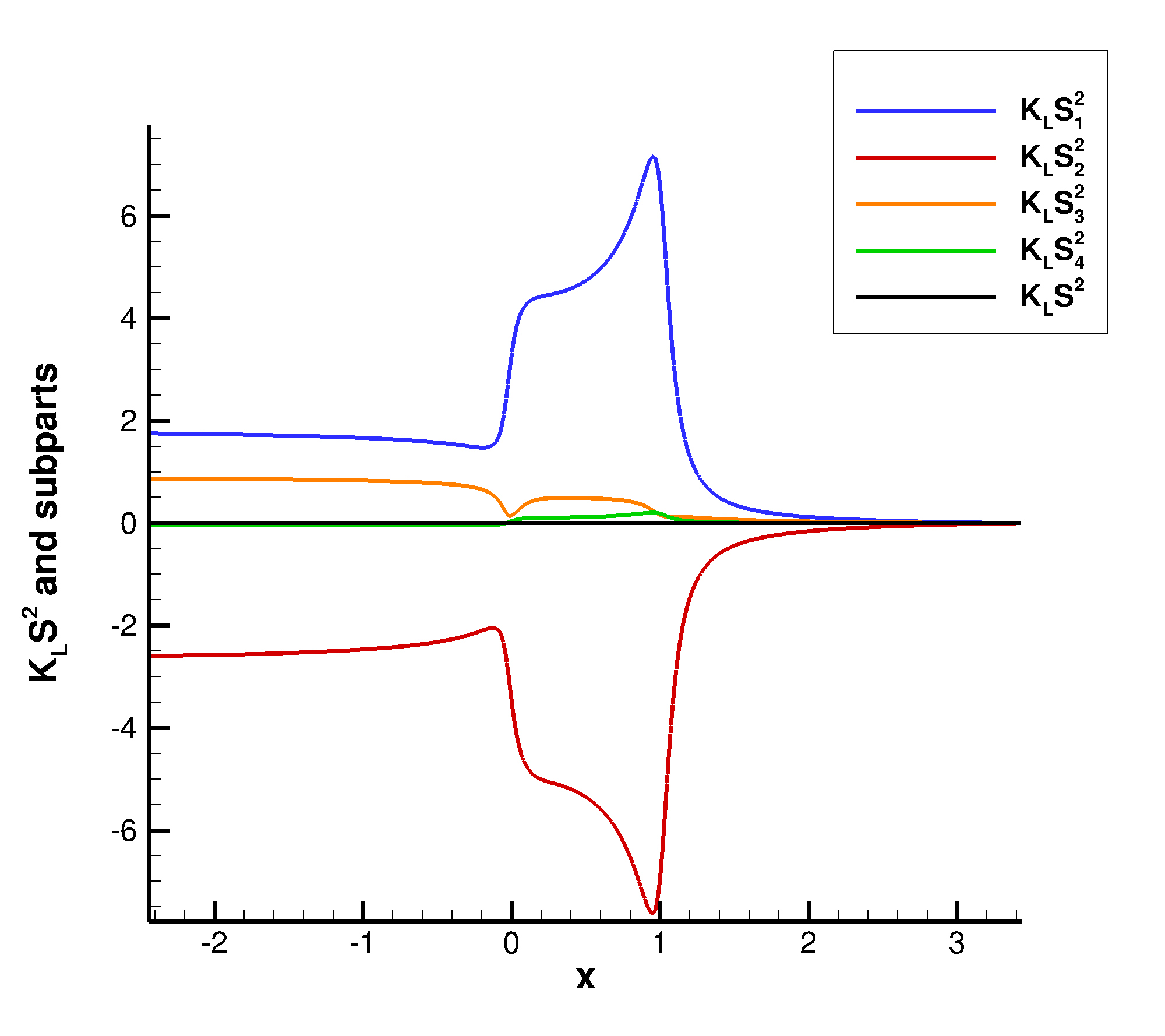}
	  \caption{$M_\infty=.40$, $\alpha=5^o$, (4097$\times$4097 mesh \cite{VasJam_10}) Numerical assessment of equation (\ref{e_atraja}) (left)
      and (\ref{e_atrajb}) (right) for the lift.
      Method of verification:  black curve should ideally coincide with the $x$ axis }
  \label{f:stream_sub}
  \end{center}
\end {figure}

We expect relations (\ref{e_atraja}) and (\ref{e_atrajb}) to be valid along
\colb{the trajectories of a subsonic flow}. The retained
flow conditions have been $M_\infty=0.4~,~ \alpha=5^o$. The  $K_DS^1$,
$K_DS^2$,  $K_LS^1$,  $K_LS^2$ integrals and their subparts are
calculated along the trajectory passing through $(c,0.1c)$. The integration indeed leads
to very small values of $K_DS^1$,
$K_DS^2$,  $K_LS^1$,  $K_LS^2$ along the curve w.r.t. their subparts. The results are equivalently
 good for lift and drag and are presented in figure \ref{f:stream_sub} for the lift.
 \colt{The lower left plot of figure \ref{f:str_mcurves} presents the typical aspect of subsonic lift or drag adjoints  indicating that an actuation able to significantly alter these
 QoIs is to be applied in the immediate vicinity of the wall. This property translates in weakly varying curves far from the profile in our verification plots for this test case.} 

\section{Conclusion}
%
\colt{ Ordinary Differential Equations have been derived for the adjoint Euler equations
 using in a first step the method of characteristics in 2D. The differential equations 
  satisfied along the streamtraces in 2D
   have then been extended to 3D 
   and the combination of equations method used for the derivation in this case also provides
    a simpler proof of the corresponding 2D equations. 
 All these ODEs are non linear differential equations that
 cannot be integrated for non-constant flows.}\\
The adjoint vector expresses the sensitivity of its corresponding
 scalar QoI to a local perturbation in the flow equations. Its variations in the fluid domain are often difficult to analyse as it precisely avoids the calculation of the flow perturbation
 that causes the change in the QoI.
  Nevertheless lift and drag 
   adjoint fields have been examined since a long time, \colt{and the presented equations
  for 2D problems clarify their highest values and the strong sensitivity of their associated QoIs to perturbations located along specific 
  lines 
  \cite{VenDar_02,TodVonBou_16,Loz_18,JPRenLab_22}.}\\
 These findings have been illustrated with flows, lift-adjoints, and drag-adjoints over the classical NACA0012 airfoil using a very fine mesh and a dual-consistent adjoint method.
 The conducted tests lead to very satisfactory results (although minor deviations
 in our
  transonic case close to the inlet of the upper side supersonic bubble).
  The demonstrated equations (\ref{e_atraja}),(\ref{e_atrajb}), and (\ref{edo_c14y})
  hence also provide a verification tool for discrete adjoint fields.
\\

\hspace{-5mm}{\bf SUPPLEMENTARY MATERIAL}\\
Supplementary material consists of:\\
-- five python scripts allowing to check the algebraic expressions of the $C^l_{mx}$ 
and $C^l_{my}$ coefficients w.r.t. their definition as determinants and one python script allowing
the numerical verification of the results of IIIB ;\\
-- six outputs of Maple scripts checking the rank of the sets of differential forms
satisfied along the streamtraces, \cp, and \cm.
\\

\hspace{-5mm}{\bf ACKNOWLEDGMENTS} \\
The authors express their warm gratitude to J.C. Vassberg and A. Jameson for allowing the co-workers of D. Destarac to  use their hierarchy of O-grids around the NACA0012 airfoil, as well as 
E. Hubert and B. Mourrain (Inria Aromath Project Team) for their helpful guidance
for the analytical developments using the computer algebra software Maple.
The authors also thank Florent Renac, Fulvio Sartor and
Martin Duguey for many fruitful discussions.
\\

This research did not receive any specific grant from funding
agencies in the public, commercial or not-for-profit sector. 
\\

\hspace{-5mm}{\bf DATA AVAILABILITY}\\
The data consists of the three flow-fields about the NACA0012 airfoil and the corresponding lift-adjoint and drag-adjoint. They are available as Tecplot binary or formatted files from the corresponding author.
\\

\hspace{-5mm}{\bf AUTHOR CONTRIBUTIONS}\\
Jacques Peter: Investigation (lead); Methodology (lead); Validation (equal); Writing  original draft (lead). Jean-Antoine D\'esid\'eri: Investigation (supporting); Methodology (supporting);  Validation (equal); Writing  original draft (supporting); Writing  review and editing.\\


\appendix

\section{\colb{Calculation and expressions of the $C_{ij}^k$ coefficients}}
%
Two or three of the vectors  $(-B_1+t A_1)$, $(-B_2+ t A_2)$, $(-B_3+t A_3)$ and $(-B_4+t A_4)$ appear in the formulas
 of the $C_{ij}^k$ coefficients expressed as the determinant of a $4 \times 4$ matrix. They may be precalculated as
 \bes
-B_1+t A_1 =  \begin{bmatrix}
0  \\
t   \\
-1   \\
0 
\end{bmatrix}
\qquad
-B_2+t A_2 =  \begin{bmatrix}
t~\gamo~E_c-u\ww  \\
-v + t (3-\gam)u   \\
-u - t \gamo v   \\
t \gamo
\end{bmatrix}
\ees
\bes
-B_3+t A_3 =  \begin{bmatrix}
-\gamo~E_c-v \ww  \\
\gamo u  + t v   \\
-(3-\gam)v + t u   \\
- \gamo
\end{bmatrix}
\qquad
-B_4+t A_4 =  \begin{bmatrix}
  (\gamo~E_c-H) \ww  \\
  t H -\gamo u \ww\\
-H - \gamo v \ww   \\
\gam \ww
\end{bmatrix}
\ees
  The expressions  of the $C^l_{1x}$, $C^l_{2x}$, $C^l_{3x}$, $C^l_{1y}$, $C^l_{2y}$, $C^l_{3y}$, $C^l_{4y}$ are gathered below (the $C^l_{4x}$ being given in \S 2.2). 
 We recall that, with our notations, the null differential form along the \st, \cm and \cp derived from the existence of, e.g. $(\partial \psi_i/\partial y)$, reads 
 $
    C_{iy}^1 d\psi_1 - C_{iy}^2 d\psi_2 + C_{iy}^3 d\psi_3 - C_{iy}^4 d\psi_4 = 0      
$
\\

\hspace{-6mm} $\bullet$ The coefficients  $C^l_{1x}$ are expressed below
\beas
C^1_{1x}&=& - dx^3  ~\ww ~( (2 t u + (t^2-1) v)(\gamo H + \gamo E_c + \gam v \ww) -(u+vt)(\gamo t H + \gamo u \ww + \gam v \ww t)) \\                                                      
C^2_{1x}&=& dx^2 dy ~ \gamo ~ \ww ~ (u^2 + v^2) ~H \\      
C^3_{1x}&=& - dx^2 dy ~\gamo ~ \ww ~t ~ (u^2 + v^2) ~H  \\      
C^4_{1x}&=& dx^2 dy  ~\gamo ~ \ww ~ (u + v t) ~H^2                                                    
\eeas
 $\bullet$ The coefficients $C^l_{2x}$  are expressed below 
\beas
C^1_{2x}&=& -dx^2 dy  ~\ww ~(\gamo H + \gamo E_c + \gam v \ww) \\                                                      
C^2_{2x}&=& -dx^3 ~\ww ~( \gamo u^2 v -u v^2 t + \gam v^3  - \gamo (u t + v) ( E_c + H)) \\      
C^3_{2x}&=& dx^2 dy ~\ww ~ ( \gamo u^2 v -u v^2 t + \gam v^3  -\gamo (u t + v) ( E_c + H)) \\      
C^4_{2x}&=& dx^2 dy ~H ~\ww ~ (\gamo H + \gamo E_c + \gam v \ww)                                        
\eeas
 $\bullet$  The coefficients $C^l_{3x}$ are expressed below
\beas
C^1_{3x}&=& dx^2 dy  ~\ww ~(t \gamo H + t \gamo E_c - \gam u \ww) \\                                                      
C^2_{3x}&=& -dx^2 dy ~\ww ~( (\gam+1) u^2 v - t \gamo u v^2 + (v + u t) ( \gamo E_c +\gamo H -\gam u^2))
  \\      
C^3_{3x}&=& dx^3  ~\ww ~ ( (t v^2 - u \ww)(-\gamo u +(\gam+1) t v) - (u t-(1+2 t^2)  v) (\gamo E_c + \gamo H -\gam v^2)) \\      
C^4_{3x}&=& -dx^2 dy ~H ~\ww ~ (t \gamo H + t \gamo E_c - \gam u \ww)                    
\eeas
 $\bullet$ The coefficients $C^l_{1y}$  are expressed below 
\beas
C^1_{1y}&=& dx^3  ~\ww ~\left( (\gamo+\gam t^2) u^3 -2 u^2 v t - \gamo t \ww H  + (\gam+\gamo t^2) u v^2 -\gamo u~ (1+t^2)~ E_c \right) \\                                                      
C^2_{1y}&=& -dx^3 \gamo ~ \ww ~ (u^2 + v^2) ~H \\      
C^3_{1y}&=&  dx^3 ~\gamo ~ t~ \ww ~ (u^2 + v^2) ~H  \\      
C^4_{1y}&=& -dx^3  ~\gamo ~ \ww ~ (u + v t) ~H^2                                                      
\eeas
 $\bullet$ The coefficients  $C^l_{2y}$  are expressed below
\beas
C^1_{2y}&=&  dx^3  ~\ww ~(\gamo H + \gamo E_c + \gam v \ww) \\                                                      
C^2_{2y}&=& -dx^3 ~\ww ~( (v\ww+u^2)(vt-(\gamo+\gam t^2)u)-(vt-(2+t^2)u))(\gamo E_c + \gamo H +\gam v\ww) ) \\      
C^3_{2y}&=& -dx^3 ~\ww ~ ( \gamo u^2 v -u v^2 t + \gam v^3  -\gamo (u t + v) ( E_c + H)) \\      
C^4_{2y}&=& -dx^3 ~H ~\ww ~ (\gamo H + \gamo E_c + \gam v \ww)       
\eeas
$\bullet$ The coefficients  $C^l_{3y}$ are expressed below   
\beas
 C_{3y}^1 &=&  -dx^3 ~\ww ~(\gamo t E_c + \gamo t H -\gam u \ww)\\
 C_{3y}^2 &=&   dx^3 ~\ww ~( (\gam+1) u^2 v - t \gamo u v^2 + (v + u t) ( \gamo E_c +\gamo H -\gam u^2)) \\
 C_{3y}^3 &=&  -dx^3 t \ww ( (\gam+1) u^2 v - t \gamo u v^2 + (v + u t) ( \gamo E_c +\gamo H -\gam u^2)) \\
 C_{3y}^4 &=&   dx^3 \ww ~H~ \left( \gamo t E_c + \gamo t H -\gam u \ww \right)
 \eeas
$\bullet $ The coefficients  $C^l_{4y}$  are expressed below 
\beas
 C_{4y}^1 &=&  -dx^3  ~ \gamma_1 ~\ww ~(u+v t)\\
 C_{4y}^2 &=&   dx^3  ~ \gamo ~ \ww ~ (u^2 + v^2) \\
 C_{4y}^3 &=&  -dx^3  ~ \gamo ~\ww~ t~ (u^2 + v^2)\\
 C_{4y}^4 &=&   dx^3 \ww  \left( (\gamo+\gam t^2) u^3 -2 u^2 v t - \gamo t \ww H  + (\gam+\gamo t^2) u v^2 -\gamo u~ (1+t^2)~ E_c \right)
 \eeas
%
\section{General properties of the $C^l_{mx}$ and $C^l_{my}$ coefficients}
%
Equation (\ref{e_cbase}) refers to the limit of small space steps and the search of
characteristic curves ; nevertheless, the expressions of the
$C^l_{mx}$ and $C^l_{my}$  coefficients may be considered for an arbitrary direction 
 and an arbitrary norm of vector $(d x, d y)$.
 Without any assumption linking $(d x, d y)$ and  $(u,v)$, the relations between the 
 coefficients  of the same differential forms are:
\bea
 C^4_{1x} &=&  H C^1_{1x} \quad C^3_{1x} = t C^2_{1x}~~;~     
 C^4_{2x} = -H C^1_{2x}  \quad C^3_{2x} = -t C^2_{2x} ~~;~    
 C^4_{3x} = -H C^1_{3x} ~~;~ C^3_{4x} = -t C^2_{4x}  \qquad  \label{e_inx} \\
 C^4_{1y} &=&  H C^1_{1y} \quad C^3_{1y} = t C^2_{1y} ~~;~ 
 C^4_{2y} = -H C^1_{2y} ~~;~~
 C^4_{3y} = -H C^1_{3y}  \quad C^3_{3y} = -t C^2_{3y} ~~;~
 C^3_{4y} = -t C^2_{4y}   \qquad  \label{e_iny} 
\eea
Besides, twelve of the sixteen coefficients of the differential forms for the $x$ and $y$ derivatives are proportional by a $(-t)$ factor: 
\bea
C^2_{1x} &=& - t ~C^2_{1y} \qquad \quad  C^3_{1x} = - t ~C^3_{1y} \qquad \quad  C^4_{1x} = - t ~C^4_{1y} \qquad \quad \label{e_cxy1} \\
C^1_{2x} &=& - t ~C^1_{2y} \qquad \quad  C^3_{2x} = - t ~C^3_{2y} \qquad \quad  C^4_{2x} = - t ~C^4_{2y} \qquad \quad \label{e_cxy2} \\
C^1_{3x} &=& - t ~C^1_{3y} \qquad \quad  C^2_{3x} = - t ~C^2_{3y} \qquad \quad  C^4_{3x} = - t ~C^4_{3y} \qquad \quad \label{e_cxy3} \\
C^1_{4x} &=& - t ~C^1_{4y} \qquad \quad  C^2_{4x} = - t ~C^2_{4y} \qquad \quad  C^3_{4x} = - t ~C^3_{4y} \qquad \quad \label{e_cxy4} 
\eea
Finally, the $C^1_{1l}$ and  $C^4_{4l}$ coefficients are equal
\beq
 C^1_{1x}= C^4_{4x}  \qquad \qquad \qquad C^1_{1y} =C^4_{4y}
\eeq
\section{Streamstrace ODEs in dimension 3}
%
From 2D equations (\ref{e_atraja}) and (\ref{e_atrajb}), we can
infer corresponding candidate equations for the adjoint vector  along the streamtraces in 3D: 
\bea
  E_c ~\frac{d\psi_1}{ds} +   H ( u~ \frac{d\psi_2}{ds} + v~ \frac{d\psi_3}{ds}  +  w~\frac{d\psi_4}{ds})  + ~H^2  \frac{d \psi_5}{ds}  = 0   \label{e_3traja} \\
  \frac{d\psi_1}{ds} + u ~ \frac{d\psi_2}{ds} +  v ~ \frac{d\psi_3}{ds} +  w ~ \frac{d\psi_4}{ds} +  E_c \frac{d\psi_5}{ds} = 0.    \label{e_3trajb}
\eea
Could these equations be possibly proven from the 3D Euler adjoint equations
\beq 
-A^T \frac{\partial \psi}{\partial x} - B^T \frac{\partial \psi}{\partial y}-
 C^T \frac{\partial \psi}{\partial z} =0,
 \label{e:adj3d}
\eeq
where the 3D transposed Jacobian of Euler fluxes read
$$
A^T = 
\begin{bmatrix} 
	0  &  (\gamma_1 E_c -u^2)  & -uv   & -uw  &(\gamma_1 E_c - H)u\\
	1  & -\gamma_1 u + 2 u        &  v    &   w  & (H-\gamma_1 u^2)\\
	0  & -\gamma_1 v         &  u    &  0   &-\gamma_1 uv \\
	0  & -\gamma_1 w         &  0    &  u   &-\gamma_1 uw \\
	0  & \gamma_1            &  0    &  0   &\gamma u\\
\end{bmatrix}
~~~~~~~~
B^T = 
\begin{bmatrix} 
	0  & -uv & (\gamma_1 E_c -v^2) &  -vw & (\gamma_1 E_c - H)v\\
	0  &  v  & -\gamma_1 u       &  0   & -\gamma_1 u v \\
	1  &  u  & -\gamma_1 v + 2v       &  w   & (H -\gamma_1 v^2)\\
	0  &  0  & -\gamma_1 w       &  v   & -\gamma_1 vw \\
	0  &  0  & \gamma_1          &  0   & \gamma v \\
\end{bmatrix}
$$
$$
C^T = 
\begin{bmatrix} 
	0  & -uw    & -vw   &  (\gamma_1 E_c -w^2)  & (\gamma_1 E_c - H)w\\
	0  &  w     &  0    &   -\gamma_1 u         & -\gamma_1 u w\\
	0  &  0     &  w    &   -\gamma_1 v         & -\gamma_1 v w \\
	1  &  u     &  v    &   -\gamma_1 w + 2w    & (H -\gamma_1 w^2) \\
	0  &  0     &  0    &    \gamma_1          &\gamma W\\
\end{bmatrix}
$$
 $s$ being the curvilinear abscissa along a trajectory in the direction opposite to the flow displacement, the differentiation w.r.t. $s$ may be expressed as
 $$ \frac{d}{ds} = -\frac{u}{\|U\|}\frac{d}{dx}
 -\frac{v}{\|U\|}\frac{d}{dy}-\frac{w}{\|U\|}\frac{d}{dz}$$
 First considering (\ref{e_3trajb}), the equation with the simpler coefficients,  this equation is satisfied in the fluid  domain if and only if 
\bea
 (u \frac{d\psi_1}{dx} + v \frac{d\psi_1}{dy} +  w \frac{d\psi_1}{dz} )~ 
 +u ( u \frac{d\psi_2}{dx} + v \frac{d\psi_2}{dy} +  w \frac{d\psi_2}{dz} )~\nonumber \\
 +v ( u \frac{d\psi_3}{dx} + v \frac{d\psi_3}{dy} +  w \frac{d\psi_3}{dz} )~
 +w ( u \frac{d\psi_4}{dx} + v \frac{d\psi_4}{dy} +  w \frac{d\psi_4}{dz} )~\nonumber\\
 +E_c ( u \frac{d\psi_5}{dx} + v \frac{d\psi_5}{dy} +  w \frac{d\psi_5}{dz} ) =0  \label{egoalb}\
\eea   
(where we have removed the norm of the velocity and the minus sign thanks to the homogeneity
of the equation). We now search if a combination of the
 lines of  (\ref{e:adj3d}) that would result in (\ref{egoalb}). For the required $\psi_1$ terms to appear, the combination
 $(-u L_2 -v L_3 -w L_4) $ ($L_j$ denoting the $j$-th line of (\ref{e:adj3d})) needs
  to be calculated:
\begin{center}
\begin{minipage}{8cm}
\begin{align*}
    & \left[ u \right. &   & (-2\gamma_1 E_c+2 u^2) &    & 2uv &  2uw   & &  \left.  (Hu-2\gamma_1 u E_c) \right] &\frac{\partial \psi}{\partial x} \\
 + & \left[ v \right. &   & 2u v  &    &  (-2\gamma_1 E_c+2 v^2) &  2vw   & &  \left. (Hv-2\gamma_1 v E_c) \right] &\frac{\partial \psi}{\partial y} \\
 + & \left[ w \right. &   & 2 u w &    & 2vw &  (-2\gamma_1 E_c+2 w^2)  & &  \left.  (Hw-2\gamma_1 w E_c) \right] &\frac{\partial \psi}{\partial z} = 0
\end{align*}
\end{minipage}
\end{center}
Subtracting the first line, calculating   $(-L_1-u L_2 -v L_3 -w L_4) $, almost fixes the expected
 coefficients for the
$\psi_1$ to $\psi_4$ derivatives:
\begin{center}
\begin{minipage}{8cm}
\begin{align*}
    & \left[ u \right. &   & (-\gamma_1 E_c+ u^2)  &    & uv &  uw   & &  \left.  -\gamma_1 u E_c \right] &\frac{\partial \psi}{\partial x} \\
 + & \left[ v \right. &   & u v  &    &  (-\gamma_1 E_c+ v^2) &  vw   & &  \left.  -\gamma_1 v E_c \right] &\frac{\partial \psi}{\partial y} \\
 + & \left[ w \right. &   & u w &    & vw &  (-\gamma_1 E_c+ w^2)  & &  \left.  -\gamma_1 w E_c \right] &\frac{\partial \psi}{\partial z} = 0
\end{align*}
\end{minipage}
\end{center}
 Finally forming   $(-L_1-u L_2 -v L_3 -w L_4 -E_c L_5)$ yields
\begin{center}
\begin{minipage}{8cm}
\begin{align*}
    & \left[ u \right. &   & u^2  &    & uv &  uw   & &  \left. uE_c \right] &\frac{\partial \psi}{\partial x} \\
 + & \left[ v \right. &   & u v  &    & v^2 &  vw   & &  \left. vE_c \right] &\frac{\partial \psi}{\partial y} \\
 + & \left[ w \right. &   & u w &    & vw &  w^2  & &  \left. wE_c \right] &\frac{\partial \psi}{\partial z} = 0
\end{align*}
\end{minipage}
\end{center}
  that is exactly equation (\ref{egoalb}).\\
In order to demonstrate  (\ref{e_3traja}), we first form the combination 
$-((2E_c-H) L_1 + u E_c L_2 +v E_c L_3 + w E_c L_4) $ of the lines of the 3D adjoint Euler 
equations. This results in 
\begin{center}
\begin{minipage}{8cm}
\begin{align*}
    & \left[ uE_c \right. &   & (u^2H-\gamma_1 H E_c) &    & uvH &  uw  & &  \left. -\gamma H u E_c + H^2 u \right] &\frac{\partial \psi}{\partial x} \\
 + & \left[ vE_c \right. &   & u v H &    & (v^2H -\gamma_1 H E_c) &  vw   & &  \left. -\gamma H v E_c + H^2 v \right] &\frac{\partial \psi}{\partial y} \\
 + & \left[ wE_c \right. &   & u w H &    & vwH &  (w^2H-\gamma_1 H E_c) & &  \left. -\gamma H w E_c + H^2 w \right] &\frac{\partial \psi}{\partial z} = 0
\end{align*}
\end{minipage}
\end{center}
Adding $- E_c~H~L_5$ yields
\begin{center}
\begin{minipage}{8cm}
\begin{align*}
    & \left[ uE_c \right. &   & u^2 H &    & uvH &  uw H  & &  \left. H^2u \right] &\frac{\partial \psi}{\partial x} \\
 + & \left[ vE_c \right. &   & u v H &    & v^2H &  vw H  & &  \left. H^2v \right] &\frac{\partial \psi}{\partial y} \\
 + & \left[ wE_c \right. &   & u w H &    & vwH &  w^2 H  & &  \left. H^2w \right] &\frac{\partial \psi}{\partial z} = 0
\end{align*}
\end{minipage}
\end{center}
that is equivalent to (\ref{e_3traja}).
%
%
\section{Demonstration of the properties presented in \S3.2 }
\colt{The differentiation of the functions $\Gamma^1_S$, $\Gamma^2_S$ and $\Gamma_{C+}$
 appearing in section \S3.2 uses the expression of the adjoint 
field where the $(\alpha-\beta)$-oriented stripe is non superimposed
 with the two other:
$$\psi(x,y) = \varphi_{\alpha-\beta} (x ~\sin(\alpha-\beta) - y ~\cos(\alpha-\beta) ) \lambda^{\alpha-\beta}_0$$}
{\footnotesize \colb{
  \beas
  \frac{d \Gamma^1_{S}}{ds} &=& E_c \frac{d\psi_1}{ds} +Hu 
  \frac{d \psi_2}{ds}+ Hv \frac{d\psi_3}{ds} + H^2 \frac{d\psi_4}{ds} \\
\frac{d \Gamma^1_{S}}{ds} &=& \left( 
 E_c( \frac{c}{\rho} ( 1 + \frac{\gamma_1}{2} M^2 )) 
+Hu (\frac{1}{\rho}( \sin(\alpha-\beta) - \gamma_1\ds\frac{u}{c})) 
+Hv (\frac{1}{\rho} ( -\cos(\alpha-\beta) - \gamma_1\ds\frac{v}{c})) 
+ H^2 (\frac{\gamma_1}{\rho c}) \right)
\frac{d \varphi_{\alpha-\beta}}{d\xi}\frac{d\xi}{ds} \\   
\frac{d \Gamma_{S}}{ds} &=& \left( \frac{\gamma_1 E_c}{\rho c} 
+\frac{c^2 u}{\gamma_1 \rho} \sin(\alpha-\beta)
-\frac{c^2 v}{\gamma_1 \rho} \cos(\alpha-\beta)
-\frac{\gamma_1 H u^2}{\rho c} -\frac{ \gamma_1  H v^2}{\rho c}
+\frac{H^2 \gamma_1}{\rho c}
\right) 
\frac{d \varphi_{\alpha-\beta}}{d\xi}\frac{d\xi}{ds}\\  
(&& \textrm{using} \quad u~\sin(\alpha-\beta)-v~\cos(\alpha-\beta)+c =0 \quad \textrm{and} \quad H=E_c + C^2/\gamma_1)\\
\frac{d \Gamma^1_{S}}{ds} &=& \left( -\frac{c^3}{\gamma_1\rho } 
+\frac{\gamma_1}{\rho c} (H-E_c)^2
\right) 
\frac{d \varphi_{\alpha-\beta}}{d\xi}\frac{d\xi}{ds}\\  
(&& \textrm{using} \quad u~\sin(\alpha-\beta)-v~\cos(\alpha-\beta)+c =0)\\
\frac{d \Gamma^1_{S}}{ds} &=& 0\\   
  \frac{d \Gamma^2_{S}}{ds} &=& \frac{d\psi_1}{ds} + u 
  \frac{d \psi_2}{ds}+ v \frac{d\psi_3}{ds} + E_c \frac{d \psi_4}{ds}  \\
\frac{d \Gamma^2_{S}}{ds} &=& \left(  
\frac{c}{\rho} ( 1 + \frac{\gamma_1}{2} M^2 )
+ u (\frac{1}{\rho} ( \sin(\alpha-\beta) - \gamma_1\ds\frac{u}{c} )) 
+ v (\frac{1}{\rho} ( -\cos(\alpha-\beta) - \gamma_1\ds\frac{v}{c})) + E_c (\frac{\gamma_1}{\rho c}) \right) 
\frac{d \varphi_{\alpha-\beta}}{d\xi}\frac{d\xi}{ds} \\   
\frac{d \Gamma^2_{S}}{ds} &=& \left( \frac{\gamma_1 E_c}{\rho c}
- \frac{\gamma_1 u^2}{\rho c} - \frac{\gamma_1 v^2}{\rho c}
+\frac{\gamma_1 E_c}{\rho c} \right) 
\frac{d \varphi_{\alpha-\beta}}{d\xi}\frac{d\xi}{ds} \\ 
(&& \textrm{using}\quad u~\sin(\alpha-\beta)-v~\cos(\alpha-\beta)+c =0 )\\
\frac{d \Gamma^2_{S}}{ds} &=& 0   
\eeas  
}}
\colt{Besides 
\footnotesize{
\beas
  \frac{d \Gamma_{C+}}{ds} &=& (u+v t^+) \frac{d\psi_1}{ds} + (u^2+ v^2) 
  \frac{d \psi_2}{ds}+ (u^2+v^2)t^+ \frac{d\psi_3}{ds} + H(u+v t^+) \frac{d\psi_4}{ds} \\
\frac{d \Gamma_{C+}}{ds} &=& ( 
 (u+v t^+)( \frac{c}{\rho} ( 1 + \frac{\gamma_1}{2} M^2 )) 
+(u^2+ v^2)(\frac{1}{\rho}( \sin(\alpha-\beta) - \gamma_1\ds\frac{u}{c}))\\ 
&&+(u^2 + v^2) t^+ (\frac{1}{\rho} ( -\cos(\alpha-\beta) - \gamma_1\ds\frac{v}{c})) 
+ H (u+v t^+) (\frac{\gamma_1}{\rho c}))
\frac{d \varphi_{\alpha-\beta}}{d\xi}\frac{d\xi}{ds} 
\eeas
}
}
\colt{
\normalsize
It is easily checked that the first and fourth terms in the bracket are 
equal. After lengthy calculations using the null eigenvalues properties and then the formulas for the difference
 of $ \cos$ and difference of $\sin$, it appears that $d \Gamma_{C+}/ds =0$. 
 For the sake of brevity, the detail of the calculations is not shown here. 
}
\normalsize
\section*{Biblography}
\bibliography{bibliotemplate}{}
\bibliographystyle{nature}

\end{document}